\documentclass[journal]{new-aiaa}
\usepackage[utf8]{inputenc}
\usepackage{textcomp}

\usepackage{graphicx}
\usepackage{amsmath}
\usepackage[version=4]{mhchem}
\usepackage{siunitx}
\usepackage{longtable,tabularx}
\usepackage{subcaption}
\usepackage{comment}
\usepackage{amsfonts}
\usepackage{xcolor}
\usepackage{algorithm, algorithmicx, algpseudocode}

\usepackage{listings}%
\usepackage{alltt}%

\newcommand{\mbf}{\mathbf}
\newcommand{\f}{\frac}
\newcommand{\bs}{\boldsymbol}
\newcommand{\R}{\mathbb{R}}
\newcommand{\mc}{\mathcal}
\newcommand{\mbb}{\mathbb}
\newcommand{\RNum}[1]{\uppercase\expandafter{\romannumeral #1\relax}}

\title{Koopman-Based Approach to Non-intrusive Projection-Based Reduced-Order Modeling with Black-Box High-Fidelity Models. Part II: Application \footnote{Submitted to the editors} }

\author{S. Ashwin Renganathan \footnote{School of Aerospace Engineering}}
\affil{Georgia Institute of Technology \\
275 Ferst Drive NW, Atlanta GA 30332}

\begin{document}

\maketitle

\begin{abstract}
A methodology for non-intrusive, projection-based non-linear model reduction originally presented by Renganathan et. al. (2018)~\cite{renganathan2018koopman} is further extended towards parametric systems with focus on application to aerospace design. Specifically, we extend the method for static systems with parametric geometry (that deforms the mesh), in addition to parametric boundary conditions. The main idea is to first perform a transformation on the governing equations such that it is lifted to a higher dimensional but linear under-determined system. This enables one to extract the system matrices easily compared to that of the original non-linear system. The under-determined system is closed with a set of model-dependent non-linear constraints upon which the model reduction is finally performed. 
The methodology is validated on the subsonic and transonic inviscid flow past the NACA0012 and the RAE2822 airfoils. We further demonstrate the utility of the approach by applying it to two common problems in aerospace design namely, derivative-free global optimization and parametric uncertainty quantification with Monte Carlo sampling. Overall, the methodology is shown to achieve accuracy upto 5\% and computational speed-up of 2-3 orders of magnitude as that of the full-order model. Comparison against another non-intrusive model reduction method revealed that the proposed approach is more robust, accurate and retains the consistency between the state variables.

\end{abstract}

\section*{Nomenclature}
{\renewcommand\arraystretch{1.0}
\noindent\begin{longtable*}{@{}l @{\quad=\quad} l@{}}
$N$  & degrees of freedom in the full model \\
$M$ & number of flow snapshots of full model (varying parameters) \\
$\mbf{N}$ & non-linear operator \\
$\mbf{L}$ & linear operator \\
$R$ & residual operator of full model \\
$\mbf{u}$ & state of the full model \\
$\tilde{\mbf{u}}$ & reduced state \\
$\mbf{U}$ & matrix with stacked $\mbf{u}$ snapshots\\
$\mbf{y}$ & observable ($=g(\mbf{u})$) \\
$\tilde{\mbf{y}}$ & reduced observable \\
$\mbf{\theta}$ & parameters of the FOM \\
$h$ & equality constraint functions \\
$S$ & number of PDEs in the full model \\
$O$ & number of terms in full model operated by a differential term \\
$\mbf{V}, \mbf{\Sigma}, \mbf{W}$ & SVD matrices \\
$\mbf{\Phi}$ \& $\mbf{\Phi}_i$ & overall and observable-wise POD matrices \\
$k$ \& $k_i$  & reduced size of all and individual observables \\
$\mc{M}$  & manifold \\
$P_\infty, \rho_\infty, a_\infty, \mu_\infty, \mc{M}, \alpha$  & free-stream parameters \\
$\rho, p, \gamma$  & density, pressure and specific-heat ratio \\
$u,v$  & $x$ and $y$ velocity components\\
$H, E$ & enthalpy and internal energy\\
$\nabla$ & gradient operator\\
$C_P, C_d, C_l$ & coefficient of pressure, drag and lift of an airfoil\\
$\mc{N}$ & Gaussian distribuion \\
$\lbrace \ell, \mu, \sigma^2 \rbrace$ & hyperparameters of the Gaussian process model\\
$\mbf{R}$ & correlation matrix of Gaussian process \\
\end{longtable*}}

\section{Introduction}
\label{s:Intro}

In the design of complex aerospace engineering systems such as aircraft, rotorcraft and spacecraft, the advent of Computer Aided Engineering (CAE) and high performance computing has significantly contributed in reducing the time and cost involved in the design process. Specifically, high-fidelity mathematical models such as Computational Fluid Dynamics (CFD), have to a good extent substituted expensive physical testing with computer simulations ~\cite{Johnson2005}. However when it comes to decision making in the design of aerospace systems, such high-fidelity models have to be queried several thousands of times, which is not practical since each query could take few days to weeks to compute even on a supercomputer. Specifically, we focus on \emph{real-time decision making} via solving \emph{many-query} problems involving \emph{global optimization} and \emph{uncertainty quantification}. 

Partial Differential Equations (PDE) based models such as CFD, typically require an iterative solution whose computational cost scales in polynomial time with the number of degrees of freedom ($N$) in the spatial discretization. Practical problems of aerospace interest, typically involve $N \sim \mathcal{O}(10^6)$ degrees of freedom, and therefore their application in the many-query context is not feasible. Therefore a \emph{surrogate model} that would trade a small amount of accuracy for a significantly large gain in computational cost is needed. Such a model would enable reliable, real-time decision making thereby offering a necessary paradigm shift in the existing process in aerospace design.

We specifically focus on Reduced Order Modeling (ROM) or Model Order Reduction (MOR) which retain the underlying physical characteristics of the high-fidelity model (henceforth referred as Full Order Model (FOM)) by projecting the actual governing equations onto a suitably chosen low-dimensional subspace, for which Proper Orthogonal Decomposition (POD) ~\cite{Lumley1998} is a common technique. ROMs are quite popular for this reason, however are limited to situations only where there is access to the discrete-form of governing equations of the FOM. In situations where the governing equations are available as a black-box (such as in commercial codes), projection-based ROM is not feasible and hence a \emph{non-intrusive} technique is necessary. In such situations, the common approach taken is to circumvent the projection step all together and use the POD to directly approximate the state whose generalized coordinates are interpolated in the time/parameter space (see ~\cite{Xiao2015, ChristopheAudouzeFlorianDeVuyst2013,C.Audouze2009,Bui-Thanh2004}). Such an approach is effective in the sense that it is purely data-driven and is more general in its applicability. However they do not guarantee that the resulting ROM still satisfies the actual governing equations and the associated consistency between the state variables which is an important characteristic of MOR. For instance, independent surrogate models for the thermodynamic variables pressure ($p$), density ($\rho$) and temperature ($T$) might not satisfy the equation of state for an ideal gas $p = \rho R T$ (where $R$ here is the specific gas constant). There has also been work done in data-driven discovery of the governing equations from black-box codes, which can then be used for model reduction~\cite{brunton2016discovering, peherstorfer2016data}; which rely on trajectory data of the state in addition to initial and boundary conditions. However, in this work we are focused on static parametric systems where there is no trajectory data and further the boundary conditions might not be explicitly available. In an earlier work, Renganathan et al (2018)~\cite{renganathan2018koopman} showed that by lifting the system to a higher dimension via the \emph{Koopman} theory~\cite{koopman1931hamiltonian}, a linear but under-determined system can be obtained which can then be closed with a set of non-linear problem-specific constraints. They further showed that with such a technique the discrete linear operator can be extracted by discretizing the linear differential terms via a method such as finite volume method, at a cost that scales linearly with the grid size, $N$. However, the method was restricted only to parameters in the boundary conditions. Here, we further extend that method to apply towards systems with parametric geometry, where the grid varies at every parameter snapshot. We show the effectiveness of the approach by generating a database of ROMs for a pre-determined set of snapshots and later interpolating between the ROMs. Therefore, this work demonstrates the methodology developed in~\cite{renganathan2018koopman} with specific applications in aerospace design by extending them to more general form of parametrizations. The baseline method for comparison is chosen to be a \emph{POD+interpolation} where as mentioned before, the generalized coordinates of the POD basis set are directly interpolated in the parameter space. The author believes that this is the only feasible non-intrusive method there exists in the literature given the same scope of the present study and hence is chosen as the baseline method. Further details are provided in section~\ref{ss:POD+Kriging}.

The rest of the paper is organized as follows. The non-intrusive MOR method is first outlined in a more generic form in section~\ref{s:Method}, followed by the compressible Euler equations in section~\ref{s:gov_eqns_exp}. The model validation is shown for subsonic and transonic inviscid flow past airfoils in section~\ref{s:validation}, to demonstrate its predictive capability. The application of the method to two many-query problems in aerospace design is discussed in section~\ref{s:Application}. The conclusion section summarizes the main findings and outlines some directions for future work. 

The following notation is followed through the rest of the paper. Scalar quantities are denoted by regular-face fonts in both upper and lower case. Vectors and matrices are represented by lower-case bold-face and upper-case bold-face fonts. The same rule applies to vectors and matrices made of several vectors and block matrices respectively. Any exception to these conventions are clarified as and when they are created in the paper.

\section{Koopman-based Non-Intrusive Reduced Order Modeling}
\label{s:Method}
The Koopman theory forms the basis of the present approach where a linear representation of the non-linear system is obtained. The finite volume method is used to discretize linear differential terms as they are well suited for unstructured computational grids and are almost the standard in commercial CFD codes. As mentioned before, this approach depends on generating a database of ROMs corresponding to a pre-determined set of parameter snapshots  which are then interpolated for new realizations of the parameters outside of the training set. The ROM interpolation draws from differential geometry in order to address the manifold-embedding of the ROM system matrices. Such an approach is essential in ensuring that the fundamental properties of the system matrices are retained post-interpolation. An outline of the methodology is provided here while the reader is referred to \cite{renganathan2018methodology, renganathan2018koopman} for more details. 

Consider a static non-linear system representing the FOM and in its discretized form

\begin{equation}
\mbf{N}(\mbf{u}) = 0
\label{e:nl_static_FOM}
\end{equation}

where $\mbf{N}$ represents a non-linear operator acting on the state variable $\mbf{u} \in \R^N$; $N$ being the degrees of freedom of the FOM. Let $g:\mbb{R}^N \rightarrow \mbb{R}^N$ be a function that operates on the state (such as $g(\mbf{u}) = $ $\mbf{u}^2, \mbf{u}^3, e^{\mbf{u}}$ etc.), then we state that

\begin{equation}
\mbf{N}(\mbf{u}) \equiv \mbf{L} \left( [g_1(\mbf{u})^\top, g_2(\mbf{u})^\top, \hdots, g_O(\mbf{u})^\top]^\top \right) 
\label{e:lifted}
\end{equation}

where, $\mbf{L}$ is a linear operator acting on the the \emph{lifted} system where the $g_i(\mbf{u})$'s replace $\mbf{u}$. We call each $g_i(\mbf{u})$ an \emph{observable} following the convention of other works on the topic (particularly~\cite{Kutz2016}) and the number of such observables $O$ in \eqref{e:lifted} is dependent on the system under consideration as is illustrated in section~\ref{s:gov_eqns_exp}. We then decompose the linear operator in the above equation as

\begin{equation}
\mbf{L}[g_1(\mbf{u})^\top, g_2(\mbf{u})^\top, \hdots, g_O(\mbf{u})^\top]^\top \approx \mbf{A}[g_1(\mbf{u})^\top, g_2(\mbf{u})^\top, \hdots, g_O(\mbf{u})^\top]^\top + \mbf{b}_a = 0
\label{e:Obs_Form}
\end{equation}

which follows the discretization of linear PDEs where, $\mbf{b_a} \in \R^{ON \times 1}$ is the vector that arises due to the discretization of boundary conditions in addition to lumping any source terms and $\mbf{A}\in \R^{N \times ON}$ is the matrix resulting from the discretization of the differential terms of a linear PDE. Overall, the parametric changes that deforms the mesh (such as geometry shape) are captured in $\mbf{A}$ whereas the rest (such as free-stream boundary conditions) are captured in $\mbf{b_a}$. Note that \eqref{e:Obs_Form} is linear but under-determined system and hence for uniqueness of the solution, constraints are added as discussed in section~\ref{ss:MOR_Cons}. Finally, we re-write \eqref{e:Obs_Form} by modifying the notation as $[g_1(\mbf{u})^\top, g_2(\mbf{u})^\top, \hdots, g_O(\mbf{u})^\top]^\top \rightarrow [\mbf{y}_1^\top,\hdots,\mbf{y}_O^\top]^\top = \mbf{y}$ and $-\mbf{b_a} \rightarrow \mbf{f}$, leading to 

\begin{equation}
\mbf{A}\mbf{y}= \mbf{f}
\label{e:transformed_FOM}
\end{equation}

\begin{figure}
\centering
\includegraphics[width = 7in]{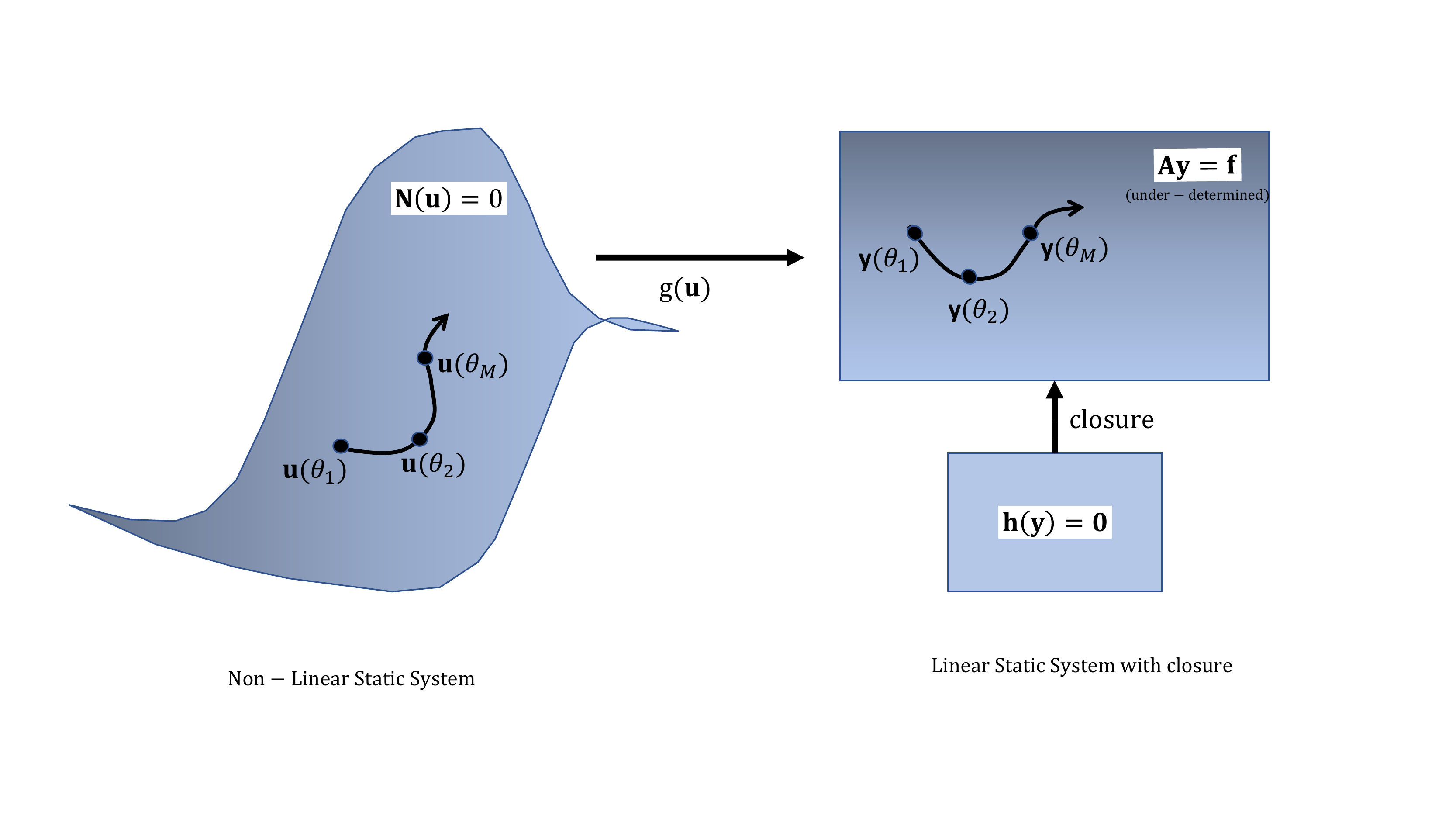}
\caption{Graphical depiction of the present methodology. The original non-linear static system is transformed to an under-determined linear system with closure~\cite{renganathan2018methodology}}
\label{f:Koopman_Method}
\end{figure}

and we work with the transformed version of the FOM in \eqref{e:transformed_FOM} to develop the ROM. Such a transformation enables us to extract $\mbf{A}$~\cite{renganathan2018koopman} and furthermore, makes the overall methodology amenable to parametric interpolation as will be illustrated in section~\ref{ss:ROM_Interp}. The overall idea behind the lifting transformation to the FOM is depicted in Figure~\ref{f:Koopman_Method} and the model reduction is performed on the transformed equations (the right hand side of the figure), which is explained in the following sub-section.

\subsection{Model Order Reduction}
\label{ss:MOR_Cons}

The total number of observables, $O$ is essentially inifinite if one were seeking a closed linear transformation of the non-linear FOM~\cite{koopman1931hamiltonian}. However we seek to find a finite $O$ which might result in an un-closed linear transformation which are then closed with a set of constraints. In this work, $O$ represents the total number of terms in the FOM that are functions of the state and are operated by a linear differential operator; each of the $g_i(\mbf{u})$'s is picked from knowledge of the FOM in its \emph{continuous} PDE form, as will be illustrated in section~\ref{s:gov_eqns_exp}. For a FOM that is a system of $S$ coupled PDEs, note  that $O \geq S$ always and $O > S$ for a non-linear system. Therefore the observables can be written

\begin{equation}
\mbf{y} = [\mbf{y}^\top_1,...,\mbf{y}^\top_S,\mbf{y}^\top_{S+1},...,\mbf{y}^\top_O]^\top
\end{equation}

Finally, to close the under-determined transformed system, we add algebraic equations that establish the non-linear consistency relationship between the observables and the state thereby providing closure. These constraints are of the form

\begin{equation}
h_i := \mbf{y}_{S+i} - f_i(\mbf{y}_1,..., \mbf{y}_S) = 0,~~i=1,\hdots,O-S
\label{e:nl_cons_generalform}
\end{equation}

where, $f_i:\R^N \rightarrow \R^N$ is a function that operates on the observables. Note that for a system of $S$ coupled PDEs, $O-S$ constraints are needed to be specified to achieve closure. Although, there is no unique way of specifying these constraints, we provide some guidelines in section~\ref{ss:MOR_Cons}. Equation~\eqref{e:nl_cons_generalform} along with \eqref{e:transformed_FOM} together form a closed system upon which model reduction is performed. To perform the projection, the truncated basis set for each observable $\mbf{y}_i$ is extracted by performing POD on the snapshot matrix of $\mbf{y}_i$ (generated by running the FOM at $M$ distinct parameter points) and are denoted $\Phi_i \in \mathbb{R}^{N \times k_i}$, i.e. 

\begin{equation}
\mbf{Y}_i = \begin{bmatrix}
\vdots & \vdots & \vdots & & \vdots \\
\mbf{y}_i^{(1)} & \mbf{y}_i^{(2)} & \mbf{y}_i^{(3)} & \hdots & \mbf{y}_i^{(M)} \\
\vdots & \vdots & \vdots & \vdots & \vdots \\
\end{bmatrix} \underset{\text{thin-svd}}{=} \mbf{V}_i \mbf{\Sigma}_i \mbf{W}_i^\top
\label{e:thin-svd}
\end{equation}

where $\mbf{y}_i^{(j)}$ is the $jth$ snapshot of observable $\mbf{y}_i$,  $\Phi_i$ is the first $k_i$ columns of $\mbf{V}_i$, $M$ is the total number of snapshots used for model training and the POD basis vectors are obtained from the thin-svd decomposition of $\mbf{Y}_i$. This leads to the trial basis matrix for the overall system defined as a block-diagonal matrix of all the $O$ POD basis set given below

\begin{equation}
    \mbf{\Phi} = \text{blkdiag} \lbrace \Phi_1, \hdots, \Phi_O \rbrace ~ \in \mathbb{R}^{ON \times k}
    \label{e:POD_basis}
\end{equation}

where  $k=k_1 + \hdots + k_O$. The \emph{reduced} observable is then given by $\tilde{\mbf{y}}_i \approx \Phi_i^\top \mbf{y}_i$. Recall that $\mbf{A}$ is non-square since it is $\in \mbb{R}^{N\times ON}$, and hence a suitable choice for the test basis for projection is $\mbf{\Psi} = \mbf{A}\mbf{\Phi}$. Note that this choice of the test basis is equivalent to a galerkin projection ($\mbf{\Psi} = \mbf{\Phi}$) on the normal equations. i.e. on $\mathbf{A}^\top \mathbf{A} \mbf{y} = \mathbf{A}^\top \mbf{f}$. Let $\mathbf{B} = \mathbf{A}^\top\mathbf{A} \in \R^{ON \times ON}$; then the projection leads to

\begin{equation}
\mbf{\Phi}^\top \mbf{B} \mbf{\Phi} \tilde{\mbf{y}} = \mbf{\Phi}^\top \mbf{A}^\top \mbf{f}
\end{equation}

Setting $\tilde{\mbf{f}} = \mbf{\Phi}^\top \mbf{A}^\top\mbf{f} \in \mathbb{R}^k$ and $\tilde{\mbf{B}} = \mbf{\Phi}^\top \mbf{B} \mbf{\Phi} \in \mathbb{R}^{k \times k}$, this leads to the reduced order model

\begin{equation}
\tilde{\mbf{B}} \tilde{\mbf{y}}=\tilde{\mbf{f}}
\label{e:ROM}
\end{equation}

 The ROM given by Equation \ref{e:ROM} is now a $k\times k$ system where $k<<N$ and is solved along with the constraints presented in Equation~\ref{e:nl_cons_generalform}, posed as a non-linear program as shown below

\begin{equation}
\begin{aligned}
& \underset{\tilde{\mbf{y}}}{\text{minimize}}
&& \frac{1}{2}\|\tilde{\mbf{B}} \tilde{\mbf{y}} - \tilde{\mbf{f}} \|_2^2 \\
& \text{s.t.}
&& \mbf{\Phi}^\top h(\mbf{\Phi} \tilde{\mbf{y}})=0
\end{aligned}
\label{e:conmin}
\end{equation}

The main hypothesis of this work is that the ROM given by Equation~\ref{e:conmin} still approximately satisfies the governing equations and this is verified in the Section~\ref{s:validation}. The optimization problem in Equation \ref{e:conmin} needs special treatment to handle the non-linear constraints which still depend on the full state of observables, and is efficiently done using the DEIM ~\cite{Chaturantabut2010}; see Appendix~\ref{A:DEIM} for details on implementation for a specific example. The ROM in \eqref{e:conmin} is solved via Sequential Quadratic Programming (SQP)~\cite{schittkowski1986nlpql} with the objective function and constraint tolerances set to $10^{-6}$ and the number of function evaluations bounded by $4\times 10^6$. The initial guess to the solution of \eqref{e:conmin} is given as the nearest flow snapshot to the test parameter at which prediction is sought.

\subsection{ROM Interpolation}
\label{ss:ROM_Interp}
\begin{figure}[htb!]
\centering
\begin{subfigure}{0.5\textwidth}
\includegraphics[width = 3in,clip = 2cm 2cm 2cm 2cm]{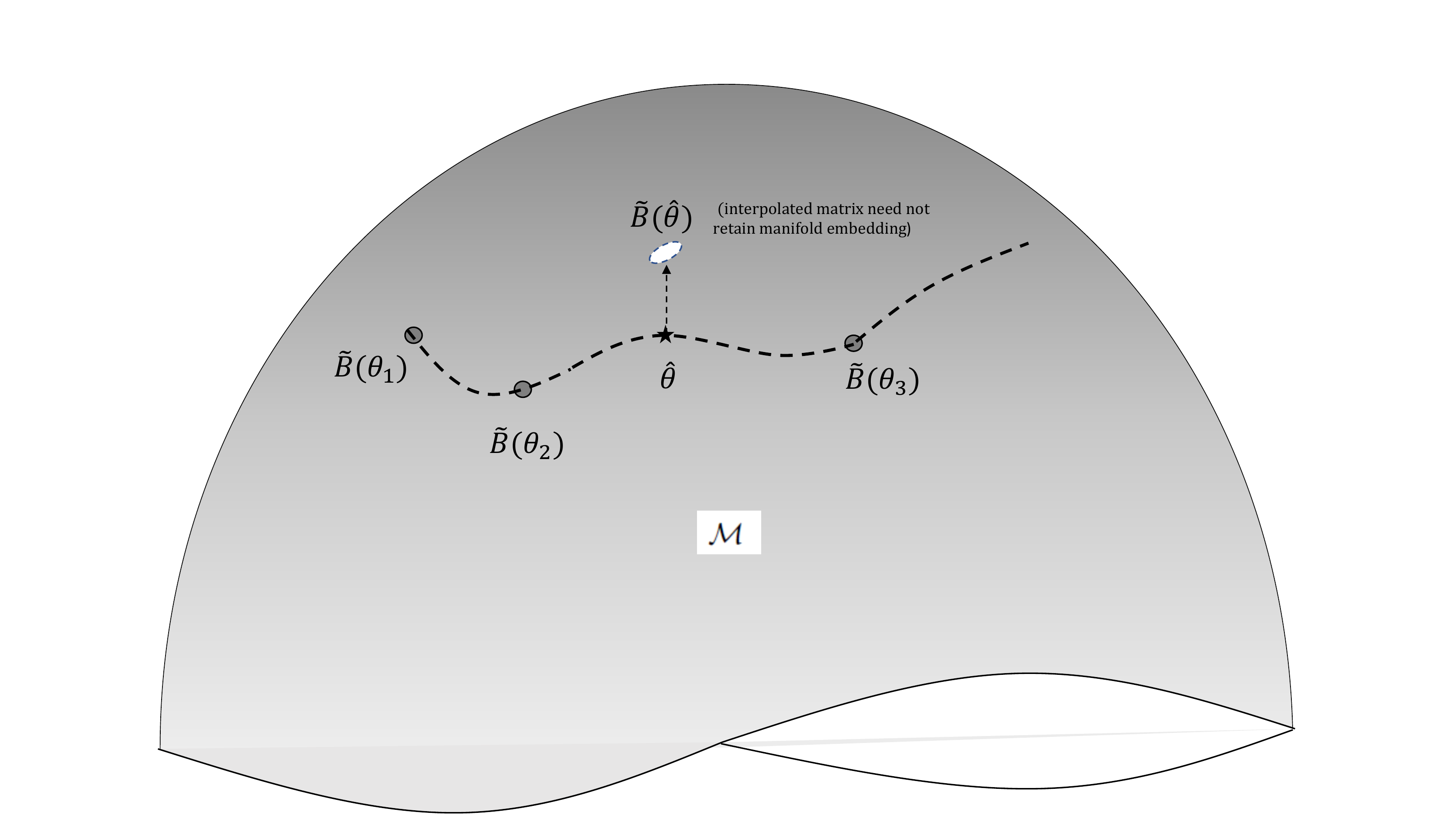}
\end{subfigure}~
\begin{subfigure}{0.5\textwidth}
\includegraphics[width = 3in,clip = 2cm 2cm 2cm 2cm]{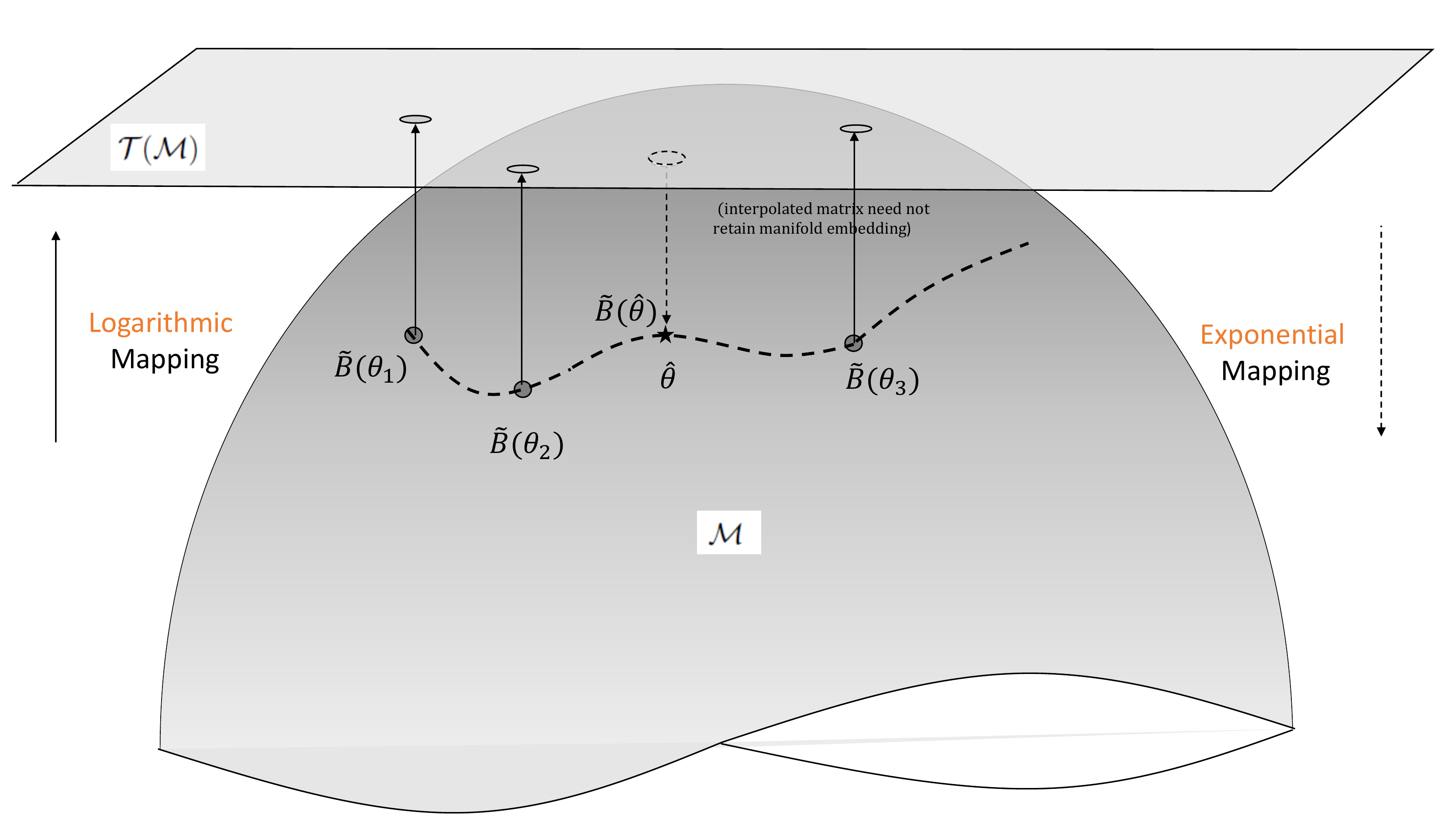}
\end{subfigure}
\caption{A graphical representation of a manifold $\mathcal{M}$ and the embedding of parametric matrices $\tilde{\mbf{B}}(\bs{\theta})$. A direct element-wise interpolation of $\tilde{\mbf{B}}$ at $\hat{\bs{\theta}}$ may not necessarily result in a matrix $\in \mathcal{M}$, and is carried out after mapping to the tangent plane, $\mathcal{T}(\mathcal{M})$}
\label{f:Manif_Intro}
\end{figure}

The proposed approach that leads to the ROM in the form of \eqref{e:ROM} corresponds to one parameter snapshot since $\tilde{\mbf{B}}$ and $\tilde{\mbf{f}}$ are parameter dependent. Therefore, the approach generates a database of ROMs for a pre-determined set of parameter snapshots, which are later interpolated to predict the state at a new parameter. The interpolation is carried out in a manner that retains the inherent structure and properties of the matrix $\mbf{\tilde{B}}$ post-interpolation. The general principle that is followed is to map the matrices to a plane that is locally tangent to the manifold in which they are originally embedded. The anchor point on the manifold, which is the point of tangency is chosen to be one of the matrices themselves. While this choice is arbitrary, in this work we use the matrix that corresponds to the nearest (in the \emph{standardized} Euclidean sense) parameter snapshot to the test parameter.
The traditional Euclidean space interpolation (where typical vector operations are valid) is then carried out in the tangent plane after which they are mapped back to the manifold. The mapping to and from the tangent plane are carried out via logarithmic and exponential relationships as depicted in Figure~\ref{f:Manif_Intro}.

The matrix $\tilde{\mbf{B}} = \mbf{\Phi}^\top \mbf{A}^\top \mbf{A} \mbf{\Phi}$ in Equation~\ref{e:conmin} is symmetric positive definite (SPD) for the following reasons. Firstly, $\mbf{A}^\top \mbf{A}$ is a co-variance matrix and hence is symmetric positive semi-definite (see~\cite{golub2012matrix}, sec. 5.3 ). Furthermore, multiplication by orthogonal matrix $\mbf{\Phi}$ of rank $k$ where $k < rank(A)$ ensures $\tilde{\mbf{B}}$ is SPD. SPD matrices of size $k \times k$ form a special group called the $SPD(k)$ ~\cite{Rahman2005, barachant2010riemannian} and the manifold they are embedded in is denoted as $\mathcal{M}$. Also, for the set of all SPD matrices $\mbf{B}~\in \mathcal{M}$, the tangent plane is the set of all \textit{symmetric} matrices, $\mbf{B'}$ ~\cite{barachant2010riemannian}. Any metric ($\mathcal{M}_f$) defined on $SPD(k)$ for any two matrices uses the following functional relationship

\begin{equation}
\mathcal{M}_f(\mbf{B}_1, \mbf{B}_2) = \mbf{B}_1^{1/2} f\left( \mbf{B}_1^{-1/2} \mbf{B}_2 \mbf{B}_1^{-1/2} \right) \mbf{B}_1^{1/2}
\label{e:SPD_f}
\end{equation}

which leads to the following results for the exponential and logarithmic mapping for $SPD(k)$~\cite{Rahman2005} where, $\tilde{\mbf{B}}_0$ is the anchor point and $\tilde{\mbf{B}'}$ is the point whose mapping is desired. The exponential mapping of $\tilde{\mbf{B}}'$ from tangent plane to $\mathcal{M}$ at $\tilde{\mbf{B}}_0 \in \mathcal{M}$, to $\mathcal{M}$ is given by 
\begin{equation}
\mbf{Exp}_{\tilde{\mbf{B}}_0} \tilde{\mbf{B}} = \tilde{\mbf{B}}_0^{1/2} \left(\tilde{\mbf{B}}_0^{-1/2} exp(\tilde{\mbf{B}}') \tilde{\mbf{B}}_0^{-1/2} \right) \tilde{\mbf{B}}_0 ^{1/2} 
\label{e:ExpSPD}
\end{equation}

and the logarithmic mapping of $\tilde{\mbf{B}} \in \mathcal{M}$ to tangent plane to $\mathcal{M}$ at $\tilde{\mbf{B}}_0 \in \mathcal{M}$: 
\begin{equation}
\mbf{Log}_{\tilde{\mbf{B}}_0} \tilde{\mbf{B}'} = \tilde{\mbf{B}}_0^{1/2} log \left(\tilde{\mbf{B}}_0^{-1/2} \tilde{\mbf{B}} \tilde{\mbf{B}}_0^{-1/2} \right) \tilde{\mbf{B}}_0^{1/2} 
\label{e:LogSPD}
\end{equation}

\begin{figure}
\centering
\includegraphics[width=.7\linewidth]{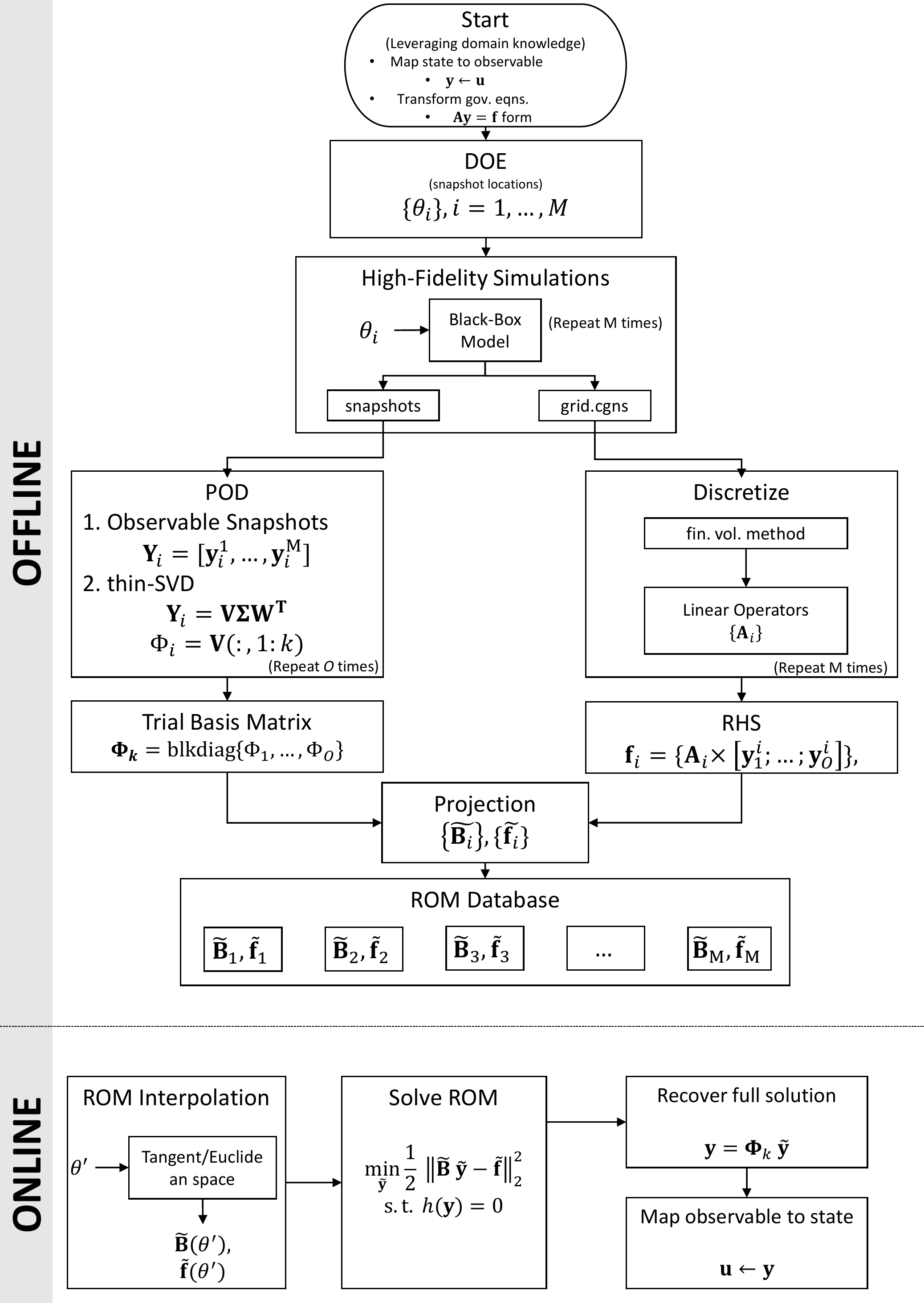}
\caption{Flowchart representation of the overall framework}
\label{f:framework}
\end{figure}

The results presented in Equations~\ref{e:ExpSPD} and \ref{e:LogSPD} are used in this work to perform mapping to and from the tangent space. Once on the tangent space the matrices are interpolated element-wise using multivariate Lagrange polynomials of $2^{nd}$ order. The overall method is summarized in Figure~\ref{f:framework}.

\clearpage

\section{Governing Equations \& Experimental Setup}
\label{s:gov_eqns_exp}
We demonstrate the methodology on the inviscid flow past airfoils, for which the NACA0012 and RAE2822 are chosen as baselines at subsonic and transonic flow regimes respectively. We begin by first illustrating the method outlined in section~\ref{s:Method} on the compressible Euler equations and then discussing the parametrization behind the chosen test cases.

\begin{table}
\caption{Free-stream conditions for the NACA and RAE test cases}
\centering
\begin{tabular}{ccc}
\hline
 & NACA & RAE \\
\hline
$p_\infty$    & 101,325 & 28,745 $Pa$\\
$\rho_\infty$ & 1.225 & 0.44 $kg/m^3$\\
$a_\infty$    & 340.296 & 301.86 $m/s$\\
$\mu_\infty$  & 1.78E-5 & 1.49E-5 $Pa-s$\\
$\mathbb{M}$  & 0.60 & 0.734 \\
$\alpha$      & 2.0 & 2.79 $deg.$ \\
\hline
\end{tabular}
\label{t:Freestream}
\end{table}

\subsection{Governing Equations}
\label{ss:FOM}
The Euler equations governing the 2D, compressible, inviscid flow past an airfoil are the governing equations on which we perform model reduction. This coupled non-linear system is solved via the commercial black-box CFD solver, STARCCM+~\cite{STARCD}. The equations in conservation form is provided in Eq.\ref{e:Euler_cons} below,

\begin{equation}
\nabla_x \mbf{F} + \nabla_y \mbf{G} = 0 
\label{e:Euler_cons}
\end{equation}

where
\begin{equation*}
\begin{aligned}
&& \mbf{F} = 
\begin{bmatrix}
\rho u \\
\rho u^2 + p \\
\rho uv \\
\rho uH
\end{bmatrix},~ \mbf{G} = \begin{bmatrix}
\rho v \\
\rho uv \\
\rho v^2 + p \\
\rho vH
\end{bmatrix} \\
&& H = E + \frac{p}{\rho}\\
&& \rho E = \frac{1}{2} \rho (u^2 + v^2) + \frac{p}{\gamma -1}
\end{aligned}
\end{equation*}

and $\nabla_x$ and $\nabla_y$ are the $x$ and $y$ components of the gradient operator $\nabla$ respectively. The following transformation is then performed  \[ [\rho u, \rho v, \rho uv, p, \rho u^2,\rho v^2, \rho uH, \rho vH]^\top  \rightarrow [y_1, y_2, y_3, y_4, y_5, y_6, y_7, y_8]^\top \] from the state variables to observables, leading to the lifted model. 
\begin{equation}
\begin{bmatrix}
\nabla_x & \nabla_y & & & & & &  \\
 &  & \nabla_y & \nabla_x & \nabla_x & & &  \\
 &  & \nabla_x & \nabla_y & &\nabla_y & &  \\
 &  & & & & & \nabla_x & \nabla_y  \\
\end{bmatrix} \begin{bmatrix}
y_1 \\
y_2 \\
y_3 \\
y_4 \\
y_5 \\
y_6 \\
y_7 \\
y_8
\end{bmatrix} = \mathbf{0}
\end{equation}

where in the above equation, empty spaces in the matrix denote zeros. The equation upon discretization leads to

\begin{equation}
\underbrace{ \begin{bmatrix}
 \mbf{G}_x & \mbf{G}_y &     &     &     &     &     &  \\
 	 &     & \mbf{G}_y & \mbf{G}_x & \mbf{G}_x &     &     &  \\
 	 &     & \mbf{G}_x & \mbf{G}_y &     & \mbf{G}_y &     &  \\
 	 &     &     &     &     &     & \mbf{G}_x & \mbf{G}_y  \\
\end{bmatrix}}_{\mbf{A}} \underbrace{ \begin{bmatrix}
\mbf{y}_1 \\
\mbf{y}_2 \\
\mbf{y}_3 \\
\mbf{y}_4 \\
\mbf{y}_5 \\
\mbf{y}_6 \\
\mbf{y}_7 \\
\mbf{y}_8
\end{bmatrix}}_\mbf{y}  = \underbrace{- \begin{bmatrix}
\mbf{b}_{a1} \\
\mbf{b}_{a2} \\
\mbf{b}_{a3} \\
\mbf{b}_{a4} \\
\mbf{b}_{a5} \\
\mbf{b}_{a6} \\
\mbf{b}_{a7} \\
\mbf{b}_{a8}
\end{bmatrix}}_{\mbf{f}}
\label{e:Euler_FOM}
\end{equation}

where, $\mbf{G}_x$ and $\mbf{G}_y$ represents the discrete version of the gradient operators $\nabla_x$ and $\nabla_y$ respectively and the empty spaces denote block matrices of zeros. The parameter-dependent $\mbf{A}$ matrix is obtained directly by discretizing the linear terms $\mbf{G}_x$ and $\mbf{G}_y$ via finite volume method. The grid is exported in the CFD General Notation System (CGNS)~\cite{CGNS} for this purpose, as mentioned in Figure~\ref{f:framework}. The snapshots $\mbf{y}$ are applied to the $\mbf{A}$ matrix and the RHS $\mbf{f}$ is extracted for each parameter. With the FOM reduced to the $\mbf{A}\mbf{y}=\mbf{f}$ form and $\mbf{A} \in \R^{4N\times 8N},~\mbf{y}, \mbf{f} \in \R^{8N}$, Equation~\ref{e:Euler_FOM} represents an under-determined system. Therefore they are closed using non-linear constraints given by Equation~\ref{e:Euler_Cons}. Notice that the constraints express the relationship between the first $S=4$ observables ($y_1$ through $y_4$) and the remaining $O-S;~(O=8)$ observables ($y_5$ through $y_8$). As mentioned in the previous section, the choice of the first $S$ observables and hence the $O-S$ constraints is non-unique. However, from experience trying out different choices in this work, it is found that the following heuristics ensure a stable transformation from the observables back to the state: (i) the terms starting from the lowest order are picked as the first $S$ observables ($\rho u, \rho v, \rho uv, p$ in this case) and (ii) one of the first $S$ observables is set to be a state variable ($y_4 = p$ in this case). It should be noted that all the observables that are in excess of the number of equations can be expressed as some function of the rest because the number of independent observables are only as many as the number of PDEs in the FOM ($S=4$). The constraints are expressed in terms of the continuous form of the state and observable as follows

\begin{equation}
\begin{aligned}
& h_1 = \rho u^2 - \frac{(\rho u)(\rho uv)}{\rho v}  \equiv y_5 - \frac{y_1 y_3}{y_2} = 0 \\
& h_2 = \rho v^2 - \frac{(\rho v)(\rho uv)}{\rho u}  \equiv y_6 - \frac{y_2 y_3}{y_1} = 0 \\
& h_3 = \rho uH - \rho u \left(E + \frac{p}{\rho} \right)  \equiv y_7 - y_1 \left(E + \frac{y_4y_3}{y_1y_2} \right)=0\\
& h_4 = \rho vH - \rho v \left(E + \frac{p}{\rho} \right) \equiv y_8 - y_2 \left(E + \frac{y_4y_3}{y_1y_2} \right)=0
\end{aligned}
\label{e:Euler_Cons}
\end{equation}

With discrete observables ($\mbf{y}_i$) all the operations in \eqref{e:Euler_cons} are performed element-wise. In all the results presented in the following sections, the POD modes constituting $99.99\%$ of the energy (cumulative fraction sum of the singular values in $\mbf{\Sigma}_i$) are retained. The error metrics used for all the results are defined as follows

\begin{equation}
\begin{split}
C_P~\text{Error} &= \frac{\left\| C_P^{FOM} - C_P^{ROM} \right \|_{\infty}}{\left\| C_P^{FOM} \right \|_{\infty}}  \times 100\\
C_d~\text{Error} &= \frac{\left| C_d^{FOM} - C_d^{ROM} \right|}{C_d^{FOM}} \times 100 \\
C_l~\text{Error} &= \frac{\left| C_l^{FOM} - C_l^{ROM} \right|}{C_l^{FOM}} \times 100
\end{split}
\label{e:error_metrics}
\end{equation}

\subsection{Test cases and parametrization}
\label{s:test_cases}
The baseline shapes are parameterized using \textit{Class Shape Transformation (CST)}~\cite{kulfan2006fundamental, kulfan2008universal}. The CST model of parametrization defines a \textit{class} function $c$ and a \textit{shape} function $s$ and the curve being parameterized is specified as their product. The main idea is that the class function serves to define a general class of geometry such as airfoils, missiles or sears-haack body, while the shape function serves to define the unique shape within a particular class of shapes (such as a NACA0012 vs RAE2822 airfoil). The class function, $c(\psi)$ is more generally defined as

\begin{equation}
c_{n_1}^{n_2} (\psi) := \psi ^{n_1}(1 - \psi)^{n_2}
\label{e:class_fn}
\end{equation}

where the variable $\psi$ represents the non-dimensional chord-wise distance. $n_1$ and $n_2$ define the specific class; for instance $n_1 = 0.5,~n_2 = 1$ and hence $\sqrt{\psi} (1 - \psi)$ defines airfoils with rounded leading edge and a sharp trailing edge~\cite{kulfan2006fundamental}. The unique shape of an airfoil is driven by the shape function, specified as follows

\begin{equation}
s(\psi) = \sum_{i=0}^{n} A_i \psi^i 
\label{e:shape_fn}
\end{equation}

where $A_i$ are the coefficients which are also the shape parameters. The NACA0012 and RAE2822 are parameterized using 6 and 8 variables respectively, whose values are given by $A_{NACA0012}$ and $A_{RAE2822}$ where the top and bottom rows correspond to the upper and lower surfaces of the airfoils. Further details of the parameterization are provided in Appendix~\ref{A:Shape}.

\begin{equation}
\begin{aligned}
A_{NACA0012} &= \begin{bmatrix}
 0.1689 & 0.2699 & 0.1387 \\
-0.1689 &-0.2699 &-0.1387
\end{bmatrix}\\
A_{RAE2822} &= \begin{bmatrix}
 0.1268 & 0.4670 & 0.5834 & 0.2103 \\
-0.1268 &-0.5425 &-0.5096 & 0.0581
\end{bmatrix}
\end{aligned}
\label{e:shape_params}
\end{equation}

\begin{figure}[htb!]
\centering
\begin{subfigure}{.5\textwidth}
\centering
\includegraphics[width = 3in]{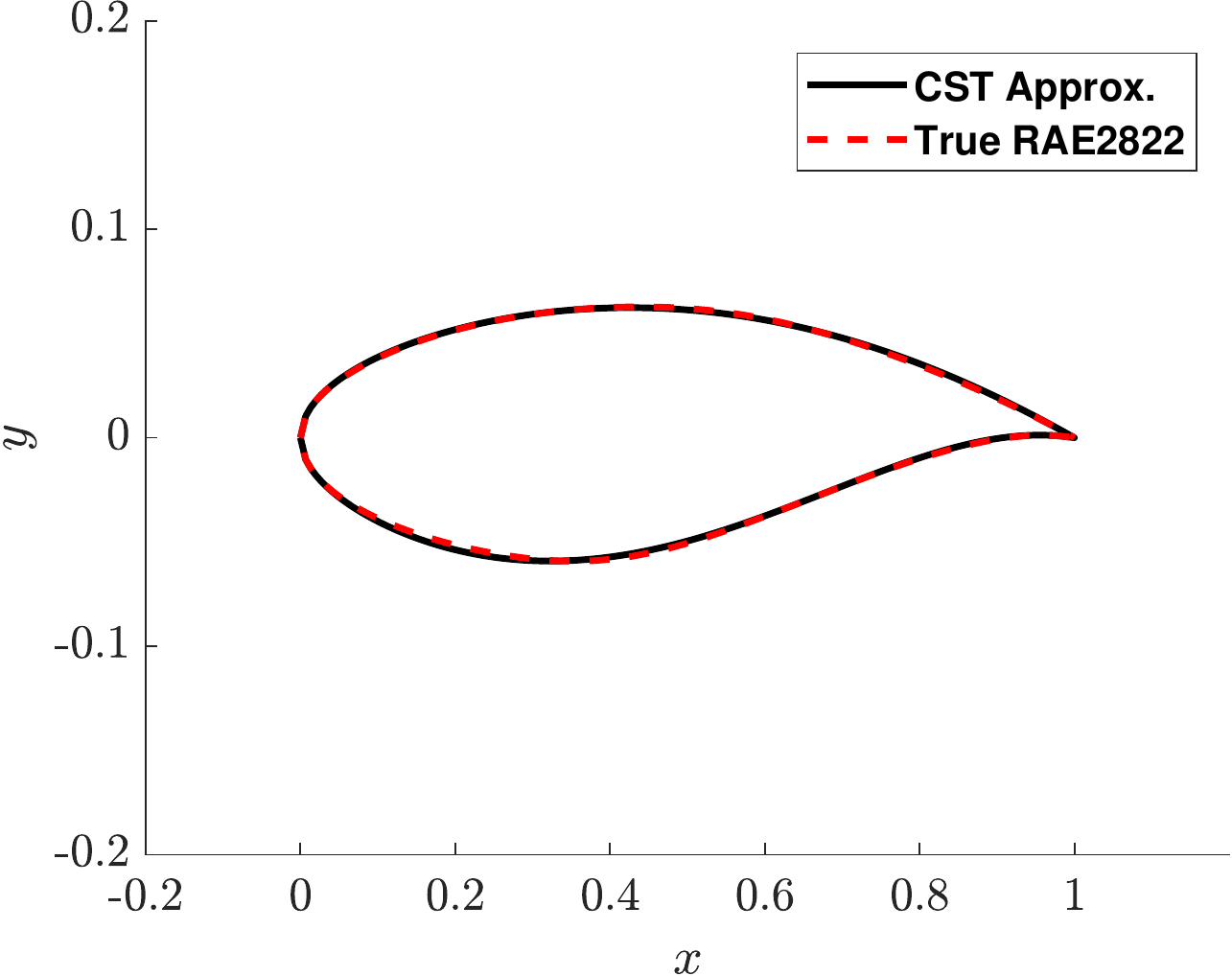}
\caption{RAE2822}
\label{f:RAE2822_CST}
\end{subfigure}~
\begin{subfigure}{.5\textwidth}
\centering
\includegraphics[width = 3in]{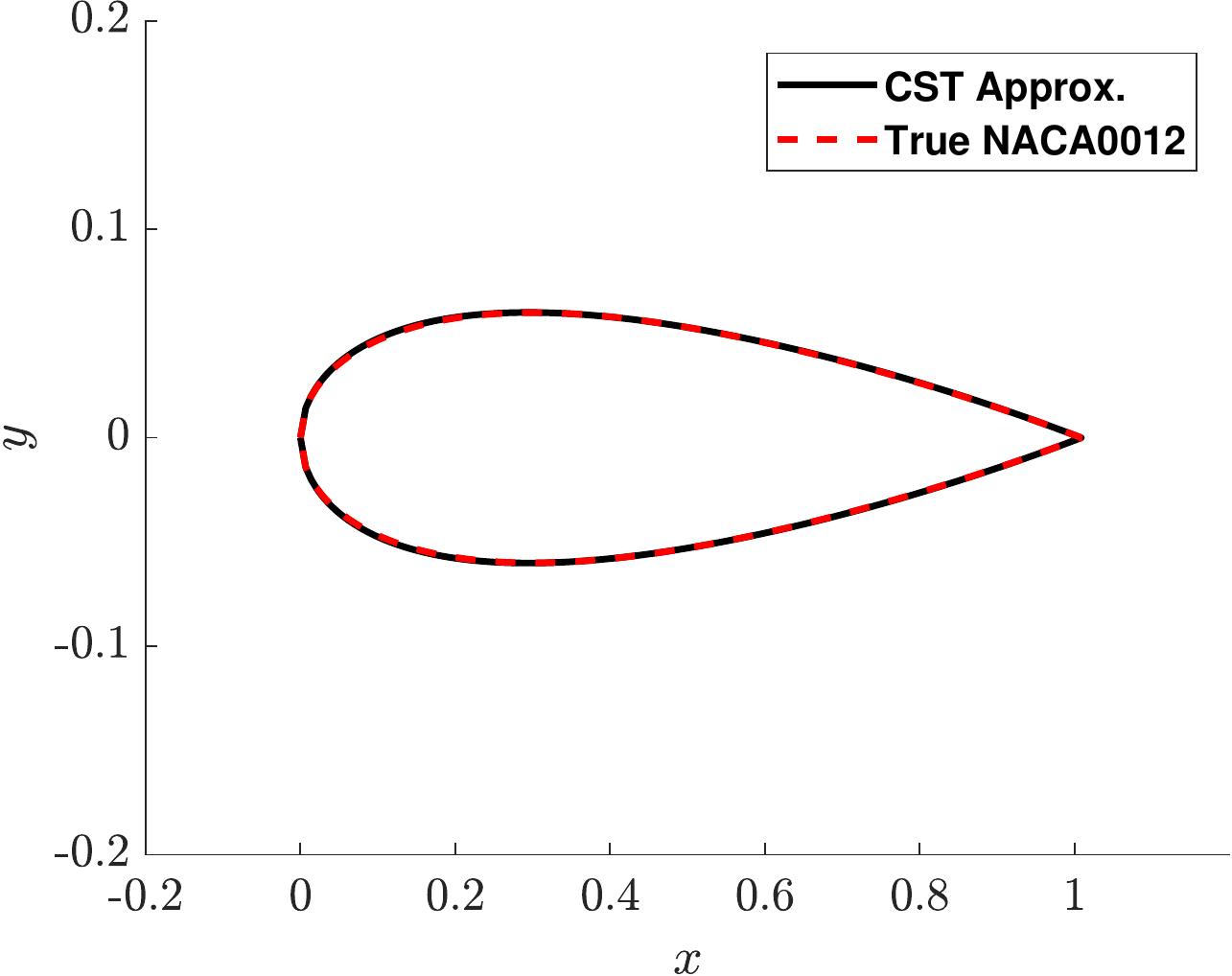}
\caption{NACA0012}
\label{f:NACA0012_CST}
\end{subfigure}
\caption{Comparison of the CST approximation against the true curve}
\label{f:CST}
\end{figure}

\section{Model Validation}
\label{s:validation}
\subsection*{NACA0012}
\begin{figure}[htb!]
\centering
\includegraphics[width=0.75\textwidth]{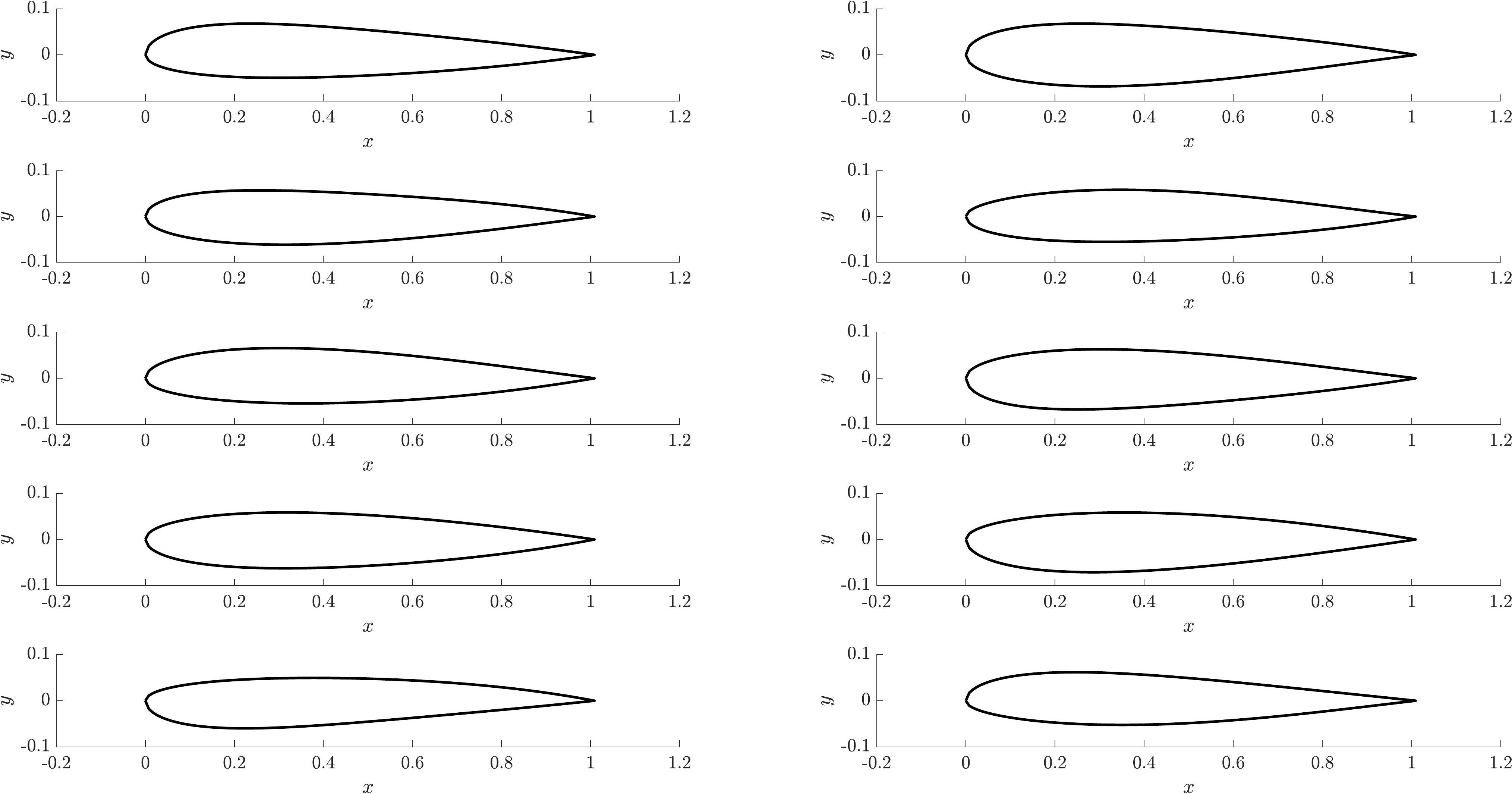}
\caption{Family of airfoils generated by perturbing (by $\pm 30\%$) the CST coefficients of the NACA0012 baseline}
\label{f:NACA0012_CST_Family}
\end{figure}

The CST coefficients representing the NACA0012 baseline is perturbed $\pm$ 30\% to generate new airfoil shapes, a sample of which is shown in Figure~\ref{f:NACA0012_CST_Family}. A total of 170 such points were generated using a Latin Hypercube design, $M = $160 of which was used in model building while the remaining was used to validate the model. The system matrix $\tilde{\mbf{B}}$ is interpolated in the tangent space to the manifold of symmetric positive definite matrices, as explained in Section~\ref{ss:ROM_Interp}, while the RHS $\tilde{\mbf{f}}$ is interpolated in the Euclidean space. In both cases a multi-variate polynomial in the Lagrange form is used for interpolation as mentioned previously (see Renganathan(2018)~\cite{renganathan2018methodology}, Algorithm 1 for details).

\begin{figure}[htb!]
\centering
\begin{tabular}{cc}
\subcaptionbox{Validation Case-1\label{1}}{\includegraphics[width = 3.25in]{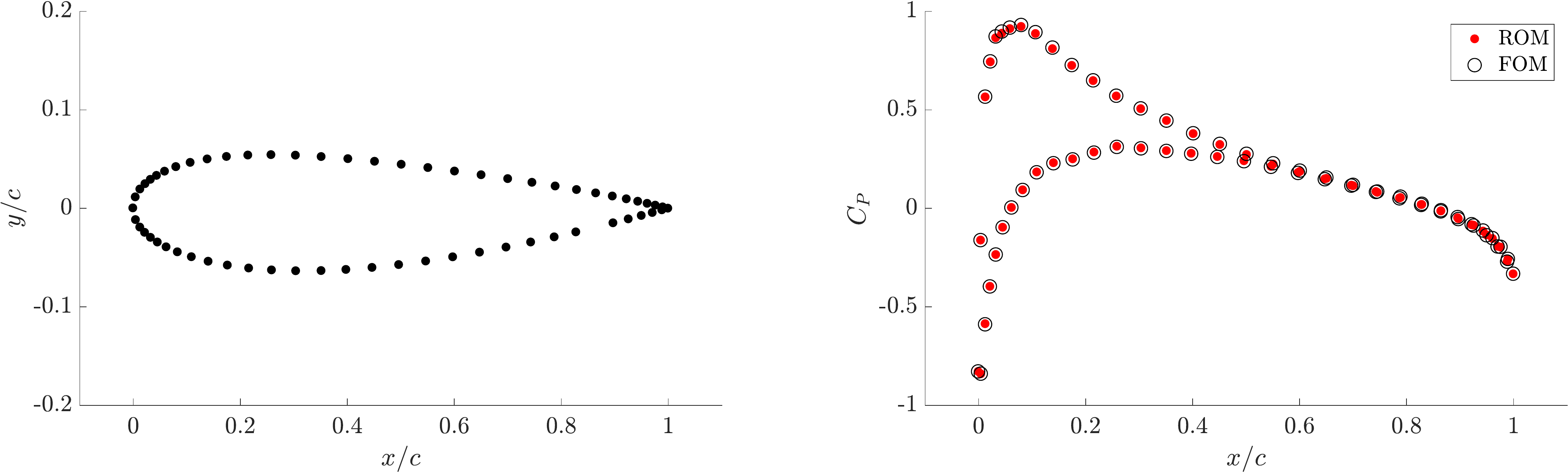}}\subcaptionbox{Validation Case-4\label{1}}{\includegraphics[width = 3.25in]{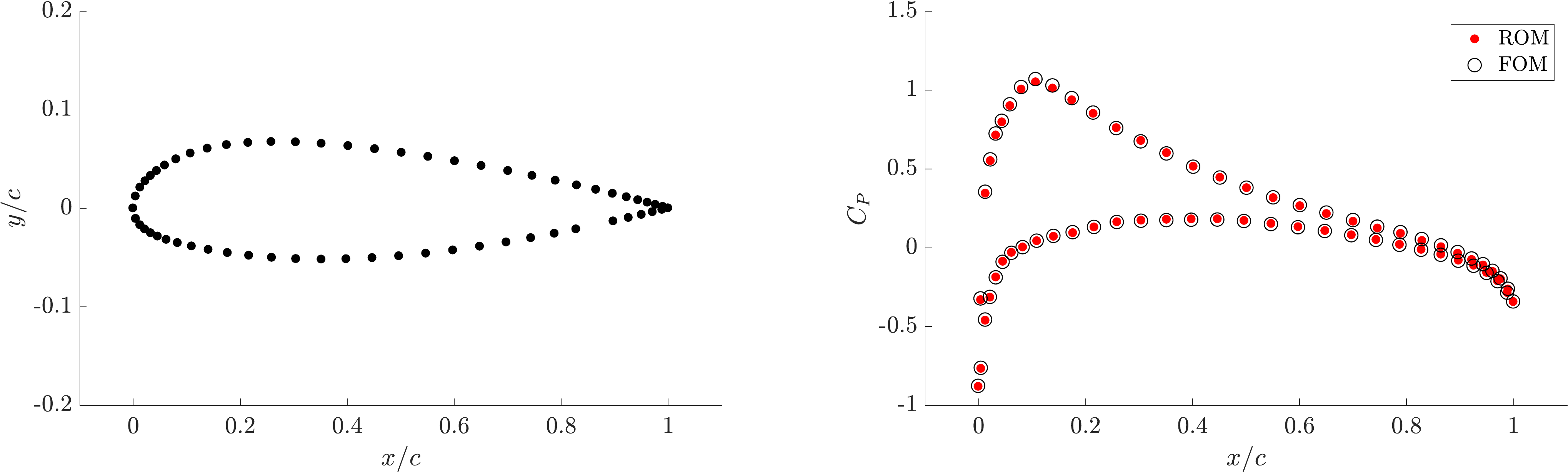}}\\
\subcaptionbox{Validation Case-1\label{1}}{\includegraphics[width = 3.25in]{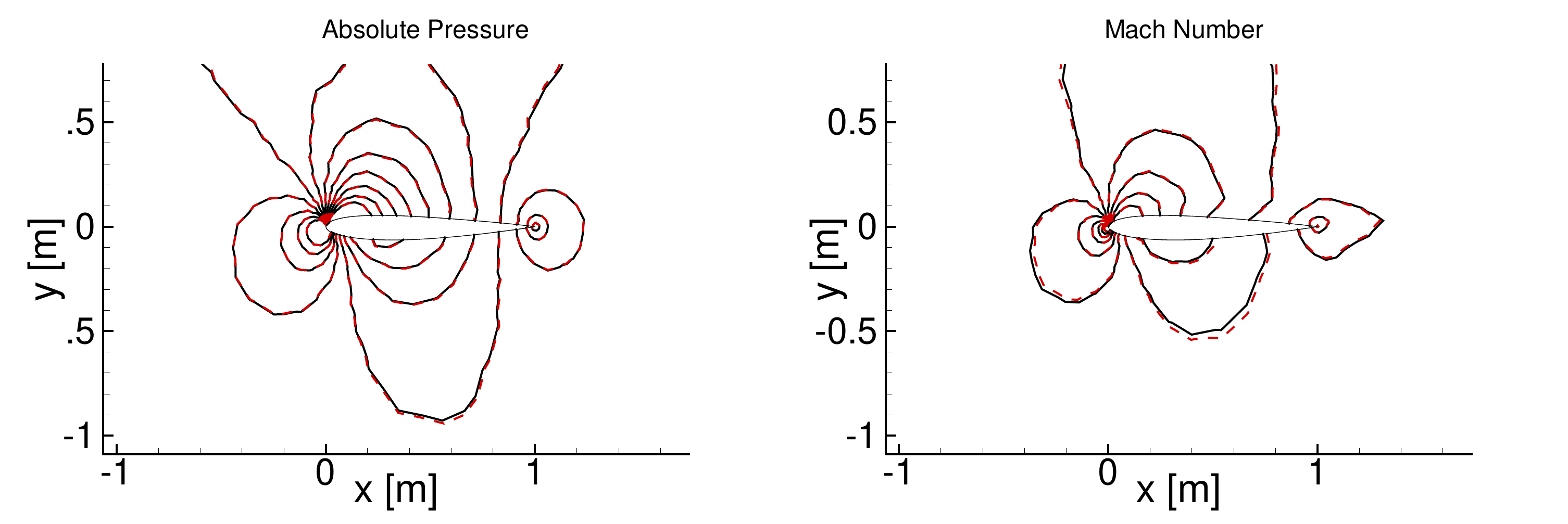}}\subcaptionbox{Validation Case-4\label{1}}{\includegraphics[width = 3.25in]{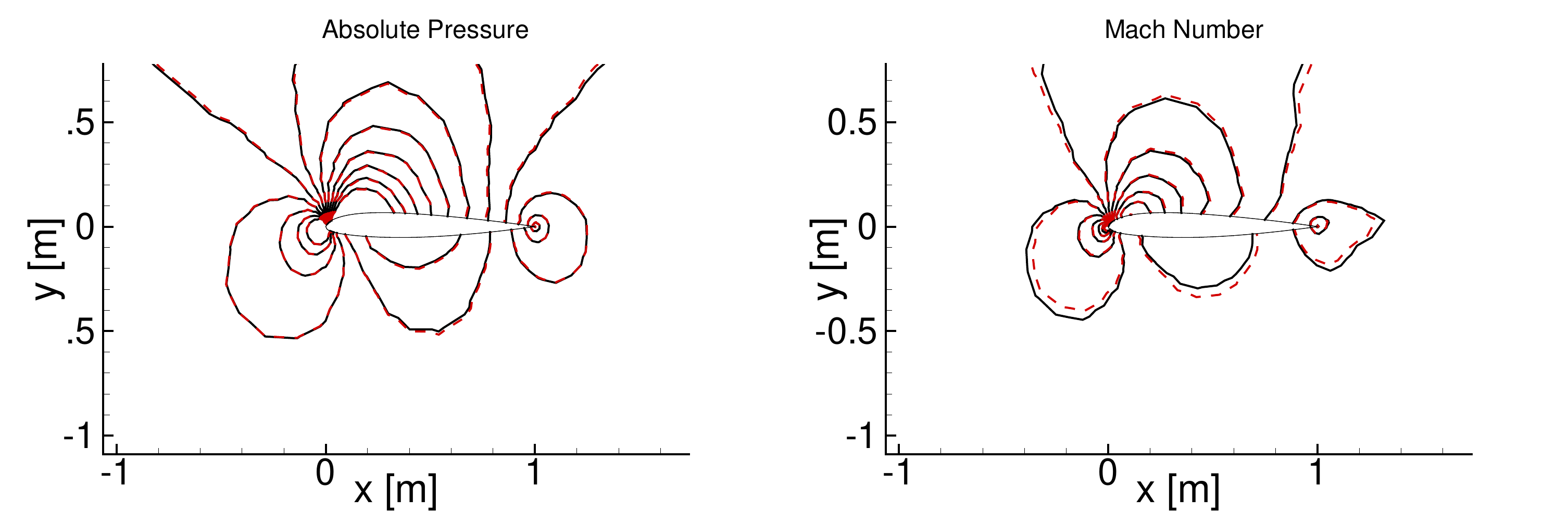}}\\
\end{tabular}
\caption{Validation of the subsonic NACA0012 test case. Top row shows airfoil shapes and $C_P$ distributions, bottom row shows Mach number and absolute pressure contours.}
\label{f:NACA_Shape_CP}
\end{figure}

 The comparison of the ROM predicted pressure coefficient on the airfoil surface against the FOM solution, for a select two of the validation set is shown in Figure~\ref{f:NACA_Shape_CP}. Across all the 10 validation cases, the maximum and average errors in $C_P$ are $\approx$ 5\% and $\approx$ 2\% respectively, whereas the computational speedup is $\approx$ 100x. Therefore, the effectiveness of proposed approach is well established. Similar accuracy was observed with the lift coefficient, which are tabulated in Appendix~\ref{A:Model_Val}. The overall flow-field comparison in terms of overlaid pressure and mach number contours is also shown in Figure~\ref{f:NACA_Shape_CP} which further re-iterate the low prediction errors observed.

\subsection*{RAE2822}
\begin{figure}[htb!]
\centering
\includegraphics[width=0.75\textwidth]{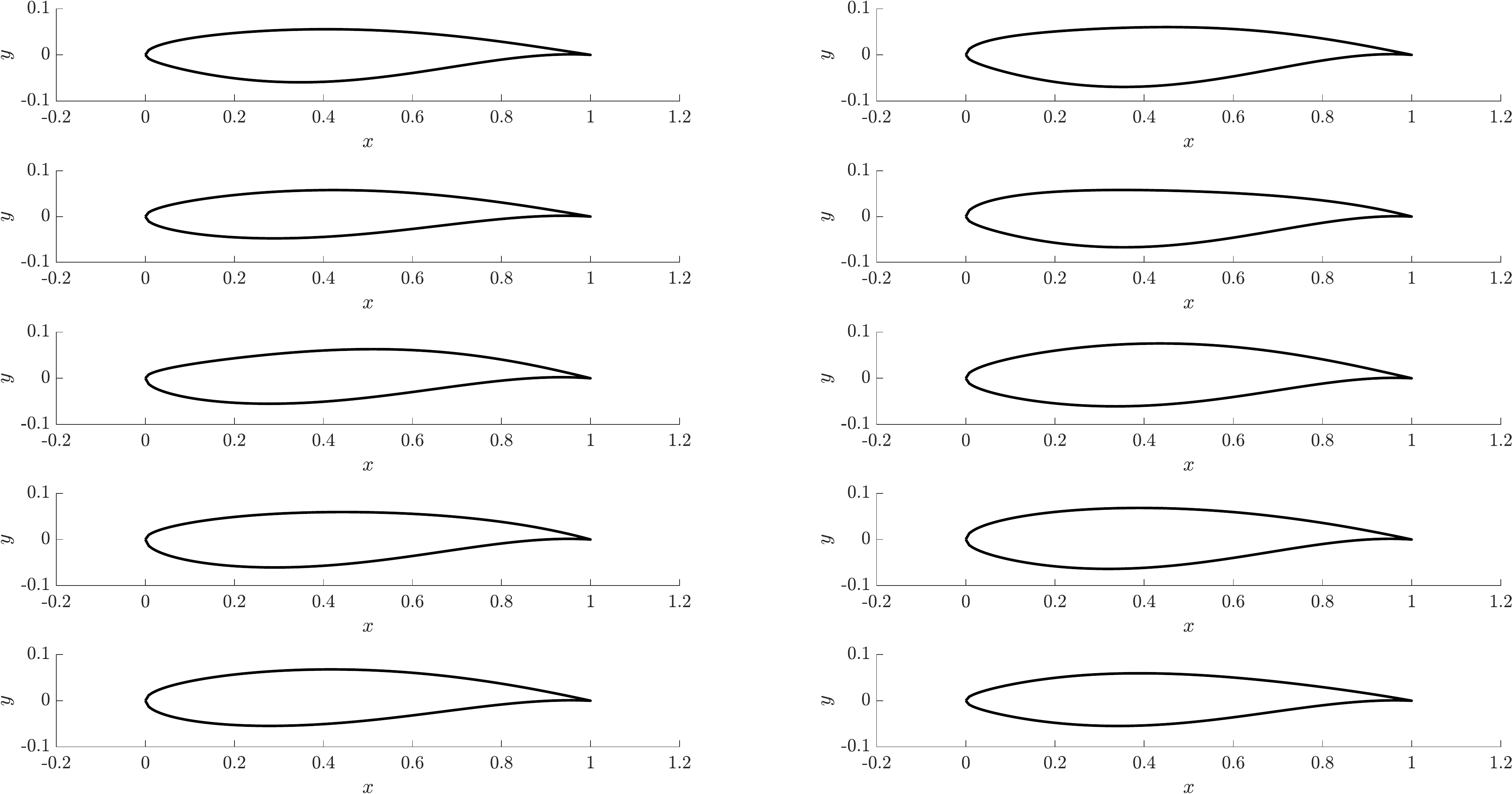}
\caption{Family of airfoils generated by perturbing (by $\pm 30\%$) the CST coefficients of the RAE2822 baseline}
\label{f:RAE2822_CST_Family}
\end{figure}

Now we demonstrate the method in the transonic regime using the RAE2822 test case. Similar to the NACA test case, CST coefficients representing the baseline are perturbed $\pm$ 30 \% to generate new airfoil shapes, a sample of which is shown in Figure~\ref{f:RAE2822_CST_Family}. A total $M = $160 snapshots were again used in model building. However, the freestream mach number for this case is set as $\mathbb{M} = 0.734$ which leads to a shock whose strength and location are affected by perturbing the shape CST coefficients. The rest of the freestream conditions are summarized in Table~\ref{t:Freestream}

\begin{figure}
\centering
\begin{tabular}{cc}
 \subcaptionbox{Val. Case-5\label{1}}{\includegraphics[width = 3.25in]{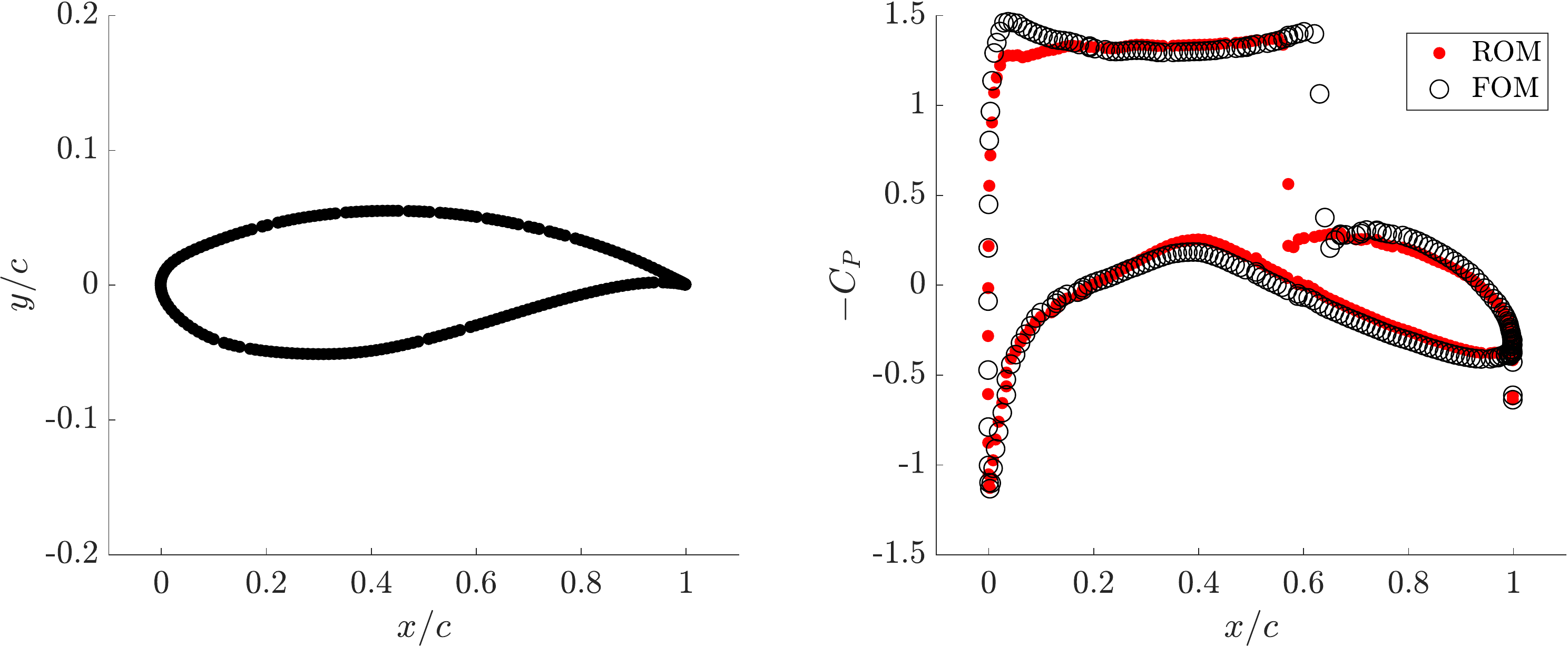}}\subcaptionbox{Val. Case-7\label{1}}{\includegraphics[width = 3.25in]{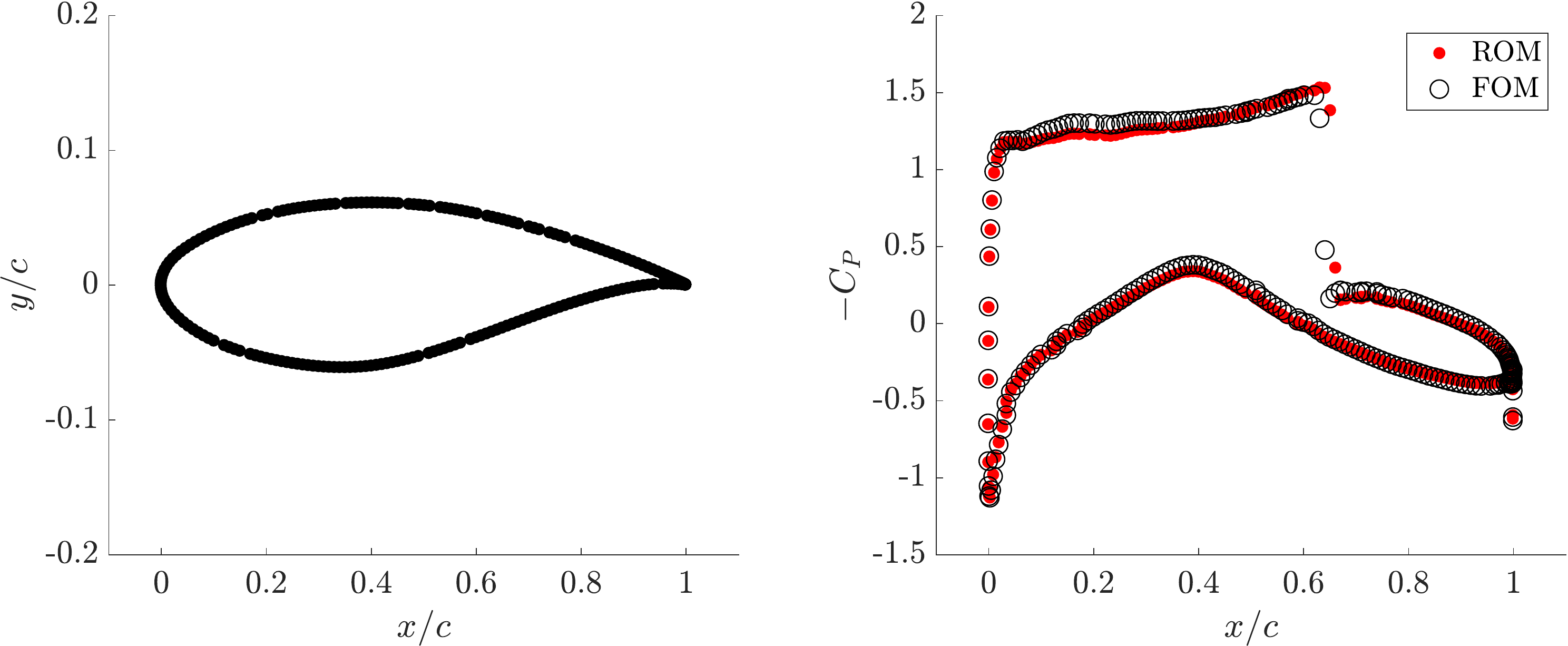}} \\
 \subcaptionbox{Val. Case-5\label{1}}{\includegraphics[width = 3.25in]{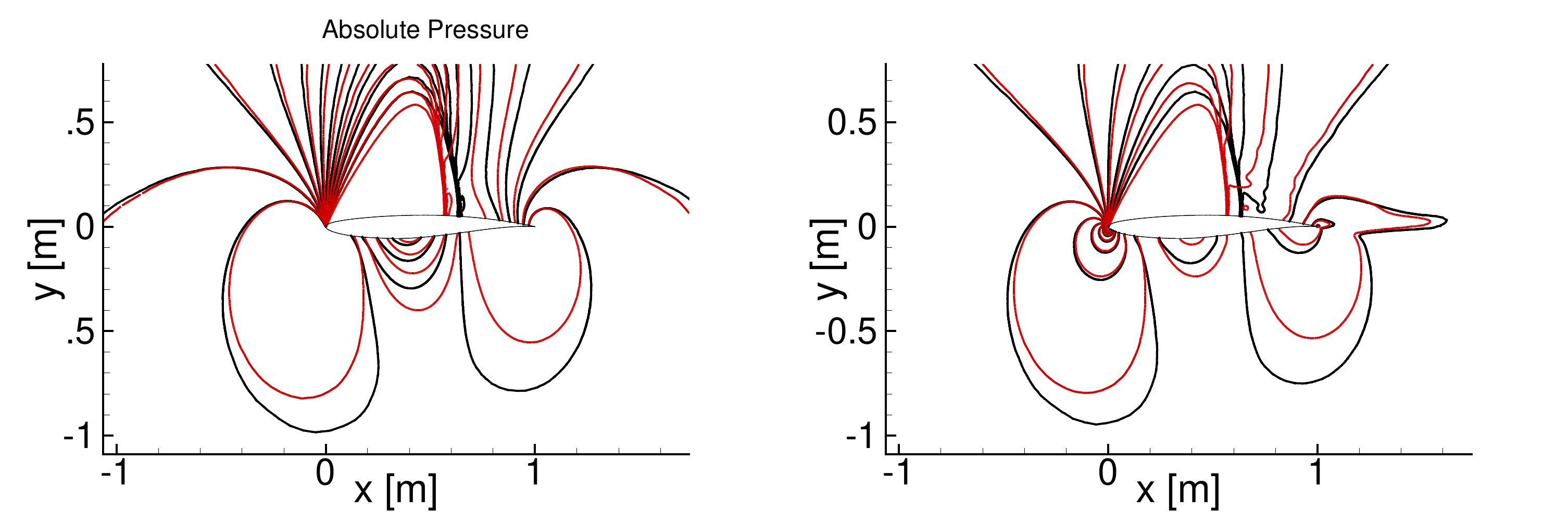}}\subcaptionbox{Val. Case-7\label{1}}{\includegraphics[width = 3.25in]{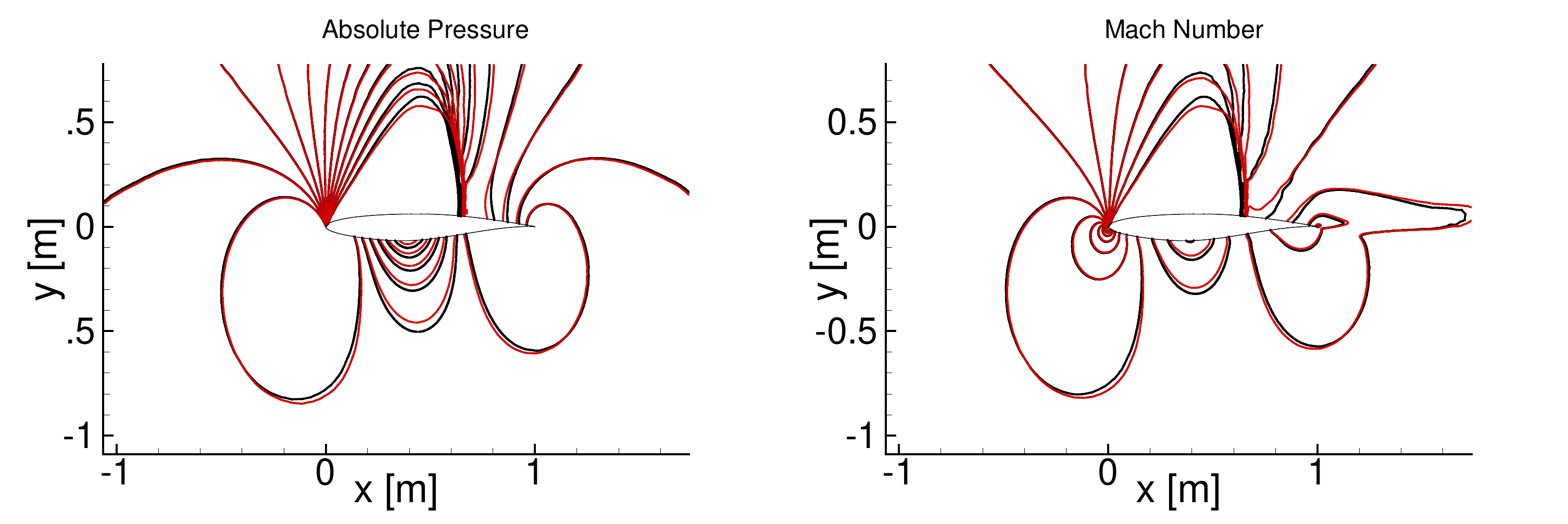}}
 \end{tabular}
\caption{Comparison of pressure coefficient $C_P$ predicted by the ROM with the true solution due to the FOM for various airfoil shapes that represent the validation cases}
\label{f:RAE_Shape_Validation}
\end{figure}

We begin by comparing the ROM predictions against FOM solutions in terms of the coefficient of pressure distributions in Figure~\ref{f:RAE_Shape_Validation}. Two specific cases are showcased in the figure to illustrate the strength and weakness of the approach in transonic regime. Overall, from all the validation cases, it is observed that the ROM does capture the shock location within a maximum error of 5\%. However, this translates in to an average $L_2$ error of 8.5\% in $C_P$ predictions
and about 16\% in the prediction of $C_d$. The average error in $C_l$ predictions were much lower at 4\%. Such a trend where the lift coefficient were predicted with greater accuracy than drag was also previously observed in \cite{renganathan2018koopman}. The main reason for this is expected to be that in the absence of viscous effects, the computation of the drag force shows greater sensitivity to discrepancy in pressure. This is because unless at high angles of attack, the projected area in the direction normal to the freestream velocity is very small leading to larger sensitivities in $C_d$ computations. However, as will be demonstrated in section~\ref{ss:UQ}, the proposed approach demonstrates more robustness in the drag prediction than competing methods.

\subsection*{Discussion}

Under subsonic shock-free conditions, the predictions are consistently under 5\% whereas under transonic conditions, the predictive capability drops. In the presence of discontinuities in the flow field such as moving shocks, a POD-based method is unable to accurately predict shock location and strength. Such a limitation has been reported by others in the literature who propose a domain-decomposition method to isolate the shock-prone regions of the flow-field where the FOM is solved whereas the ROM is used to approximate the rest of the flow field~\cite{legresley2006application, Lucia2002}. However such methods are intrusive in the sense that they might require access to source code of the FOM to effectively manipulate the flow domains which is not feasible with black-box codes. Another method reported in the literature to address this problem is to do a space-transformation of the snapshots containing moving discontinuities~\cite{welper2017interpolation}. The goal of the proposed approach are first steps towards enabling projection-based MOR for black-box static parametric systems although such enhancements can be added to it to address specific problems. Additionally, the solution to the ROM in \eqref{e:ROM} shows sensitivity to initial guess and hence can lead to a local optimum which also contributes to the performance penalty of the approach in transonic regime. Despite the limitations in the transonic regime, it has been observed that the ROM still predicts the shock within 5\% chord-length variability. Furthermore, the $C_P$, $C_d$ and $C_l$ are predicted with an average error of $14\%$, $16\%$ and $4\%$ respectively which are still useful in the conceptual and preliminary stages of the design process where some accuracy can be traded for larger gain in computational costs for faster decision making.

\section{Application: Many-Query Problems}
\label{s:Application}
We pick two problems which are pertinent to aerospace design and require several queries to the model which is not practical if one had only the FOM. The first problem is a non-convex optimization problem which greatly benefits from derivative-free optimization techniques such as the Genetic Algorithm (GA)~\cite{GoldbergGA1989} that efficiently balances exploration and exploitation to determine the global optimum. The second problem concerns the approximation of probability densities which similarly require several queries to the model.

As a baseline for comparison of the proposed approach, we use a more commonly used non-intrusive technique which replaces the projection step in MOR with a direct interpolation of the reduced state. As mentioned previously, such an approach no longer guarantees that the ROM still satisfies the governing equations of the FOM, but their accuracy has been found useful in multiple studies~\cite{ChristopheAudouzeFlorianDeVuyst2013, Xiao_RBF}. Here we replace the projection step with the Kriging~\cite{stein2012interpolation, rasmussen2003gaussian} interpolator and use that for comparison against the proposed method. In the rest of the section the proposed projection-based approach is referred as \emph{POD-Proj.} whereas the interpolation-based approach is referred as \emph{POD-Krig.}

\subsection{POD + Kriging}
\label{ss:POD+Kriging}
Recall from section~\ref{s:Method} that the full and reduced observable are related via the relationship $\tilde{\mbf{y}}_i \approx \mbf{\Phi}^\top \mbf{y}_i$. Here each element of the reduced observable $\tilde{\mbf{y}}_i$ is assumed to be a smooth function of the parameters $\bs{\theta}$ and is interpolated in the $\bs{\theta}$-space via Kriging. The Kriging approach assumes that the true response is a relaization of a Gaussian process (GP) which is completely specified by a \emph{mean} and \emph{covariance function} and estimates the hyperparameters of the GP from observations at finite number of locations. Kriging interpolates noise-free data, regresses noisy data and the prediction at unknown sites is given by a Gaussian probability distribution as opposed to a deterministic value. However it is customary to use the expected value of the predictive distribution as a surrogate for the true function.

The $jth$ element of the $ith$ reduced observable is assumed to take the form
\begin{equation}
    \tilde{y}_{ij} = \mu + z(\bs{\theta})
    \label{e:krig_underlying_fn}
\end{equation}

where $z(\bs{\theta}) \sim \mc{N}(\mbf{0}, \mbf{C})$ with $\mbf{0}$ being a vector of zeros of appropriate length, $\mbf{C}$ the covariance matrix and $\mu$ the unknown mean. The GP is assumed to be statistically \emph{stationary}~\cite{santner2003design} with homoscedastic variance which simplifies the covariance matrix as $\mbf{C}_{ij} = \sigma^2 \mbf{R}_{ij}$, where $\mbf{R}$ is the \emph{correlation} matrix. The assumption of stationarity also means that $\mbf{R}_{ij} = \mbf{R}(\| \bs{\theta}_i - \bs{\theta}_j \|)$ where $\| \cdot \|$ denotes the Euclidean distance. The elements of the correlation matrix are defined by a parametrized \emph{kernel} as $\mbf{R}_{ij} = r(\bs{\theta}_i, \bs{\theta}_j ~; \ell)$. The hyperparameters of the Kriging model are therefore $\lbrace \ell, \mu, \sigma^2 \rbrace$. In this work, we assume the stationary and isotropic \emph{squared-exponential} kernel are given by the following equation

\begin{equation}
        \mbf{R}_{ij} = r(\bs{\theta}_i , \bs{\theta}_j; \ell) = \text{exp} \left(-\f{\|\bs{\theta}_i - \bs{\theta}_j\|_2^2}{2\ell^2} \right)
    \label{e:SE_kernel}
\end{equation}

The training data for each $\tilde{y}_{ij}$ is the $jth$ row of the product $\mbf{\Sigma}_i \times \mbf{W}_i^\top$ in the thin-svd step in ~\eqref{e:thin-svd}; let that be denoted as $\tilde{\mbf{y}}^j$. Each $\tilde{\mbf{y}}^j$ is a vector of length $M$ corresponding to the parameter snapshots $\bs{\Theta} = [\bs{\theta}_1, \hdots, \bs{\theta}_M]$. Then the Kriging model assumes that the prediction at some unknown $\bs{\theta}^*$ forms a joint normal distribution with $\tilde{\mbf{y}}^j$. That is

\begin{equation}
    \begin{bmatrix}
    \tilde{\mbf{y}}^j \\
    \tilde{y}_{ij}^*
    \end{bmatrix} = 
    \sim \mathcal{N} \left(\mbf{0}, 
    \sigma^2 \begin{bmatrix}
    \mbf{R}      & r(\bs{\Theta},   \bs{\theta}^*) \\
    r(\bs{\theta}^*, \bs{\Theta}) & 1
    \end{bmatrix} \right)
    \label{e:Krig_joint_normal}
\end{equation}

Conditioning the distribution of  $\tilde{y}_{ij}^*$ in \eqref{e:Krig_joint_normal} on the observations one then obtains the posterior predictive distribution~\cite{rasmussen2003gaussian} given by

\begin{equation}
    \tilde{y}_{ij}*|\bs{\theta}^*, \bs{\Theta},\tilde{\mbf{y}}^j \sim \mathcal{N} \left(r(\bs{\theta}^*, \bs{\Theta}) \mbf{R}^{-1} \tilde{\mbf{y}}^j, \sigma^2 \left[ 1 - r(\bs{\theta}^*, \bs{\Theta}) \mbf{R}^{-1} r(\bs{\Theta},   \bs{\theta}^*) \right]\right) 
    \label{e:krig_posterior}
\end{equation}

where $r(\bs{\theta}^*, \bs{\Theta}) = r(\bs{\Theta}, \bs{\theta}^*)^\top$ is a column vector of the correlation between $\bs{\theta}^*$ and $\bs{\Theta}$ and the unknown mean $\mu = r(\bs{\theta}^*, \bs{\Theta}) \mbf{R}^{-1} \tilde{\mbf{y}}^j$ (which is also our Kriging predictor). The other two hyperparameters ($\ell, \sigma^2$) are determined by maximizing the log marginal likelihood of $\tilde{y}_{ij}*$

\begin{equation}
    \text{log}~ p\left(\tilde{\mbf{y}}^j|\bs{\Theta}, \ell, \sigma^2 \right) = -\f{1}{2}\tilde{\mbf{y}}^{j \top} \mbf{R}^{-1} \tilde{\mbf{y}}^j -\f{1}{2}\text{log}|\mbf{R}| - \frac{M}{2}\text{log}~2\pi
    \label{e:log_ML}
\end{equation}

\subsection{Inverse Design}
\label{ss:inv_design}
In aerodynamic design, a specific aerodynamic load distribution about an aerodynamic object is of interest. For instance, under incompressible flow assmptions for a finite wing, an elliptic lift distribution along the wing is known to produce the least induced drag \cite{anderson2010fundamentals}. Similarly, in the preliminary design of propellers~\cite{adkins1994design} and turbines\cite{renganathan2014validation}, a certain lift distribution along the blade is an input to the design process. In such cases the actual design (shape) parameters that produce such a load distribution is of interest. We call such a problem the \emph{inverse design}. Here we fix the desired pressure coefficient distribution as our \emph{target} and search the design space for the shape parameters that would best approximate the target. Then the predicted airfoil shape is compared against the true shape. Therefore we are interested in solving the following optimization problem

\begin{equation}
\begin{aligned} 
 \underset{\bs{\theta}} {arg\text{min}} && ~~\frac{1}{2}\| C_P(\mbf{u},\bs{\theta}) - C^{*}_P(\mbf{u},\bs{\theta})\|_2^2  \\
\text{subject to:} \\
&&  R(\mbf{u}, \bs{\theta}) = 0 \\
&& \bs{\theta}_l \le \bs{\theta} \le \bs{\theta}_u \\ 
\end{aligned} 
\label{e:Inverse_design} 
\end{equation}

where $C^{*}_P$ is the target pressure distribution and $R$ is the residual operator. Naturally we want to replace the full-order governing equations with the ROM and hence we solve the modified problem

\begin{equation}
\begin{aligned} 
 \underset{\bs{\theta}} {arg\text{min}} && ~~\frac{1}{2}\| C_P(\mbf{u},\bs{\theta}) - C^{*}_P(\mbf{u},\bs{\theta})\|_2^2  \\
\text{subject to:} \\
&& \mbf{\Psi}^\top R(\mbf{\Phi} \tilde{\mbf{u}}, \bs{\theta}) = 0 \\
&& \bs{\theta}_l \le \bs{\theta} \le \bs{\theta}_u \\ 
\end{aligned} 
\label{e:Inverse_design_ROM} 
\end{equation}

\begin{table}
\centering
\caption{Free-stream conditions for inverse design}
\begin{tabular}{cc}
\hline
\hline 
$P_{\infty}$ & 101325 Pa \\ 
$T_{\infty}$ & 288 K \\ 
$\rho_{\infty}$ & 1.225 $kg/m^3$ \\ 
$a_{\infty}$ & 340.296 m/s \\ 
$M_{\infty}$ & 0.6 \\ 
$\alpha$ & 2 $^\circ$\\
\hline 
\end{tabular}
\label{t:Inv_des}
\end{table}

\begin{figure}[htb!]
\centering
\begin{subfigure}{0.5\textwidth}
	\centering
	\includegraphics[width= 1\linewidth]{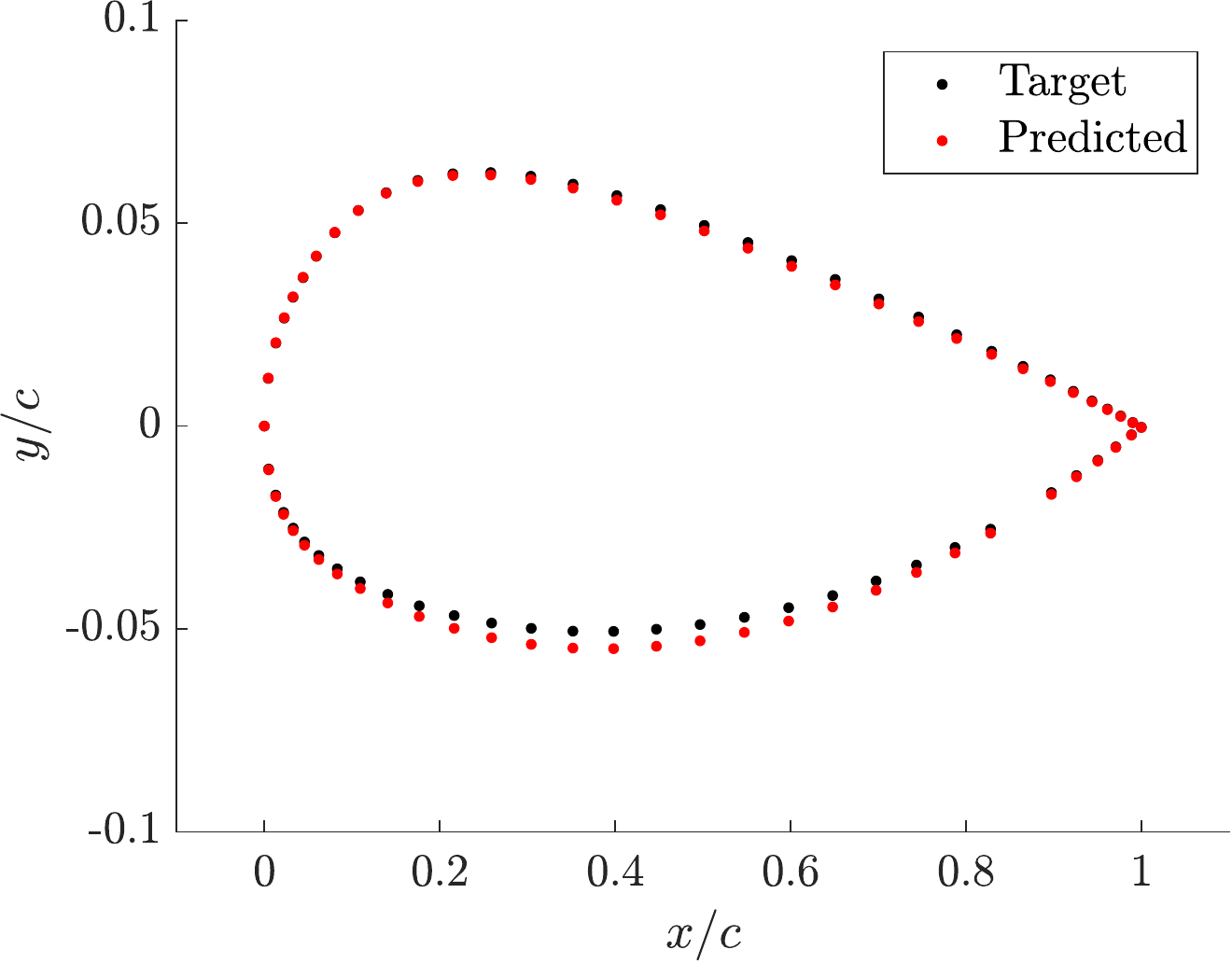}
	\caption{Airfoil shape comparison}
	\label{f:}
\end{subfigure}%
\begin{subfigure}{0.5\textwidth}
	\centering
	\includegraphics[width= 1\linewidth]{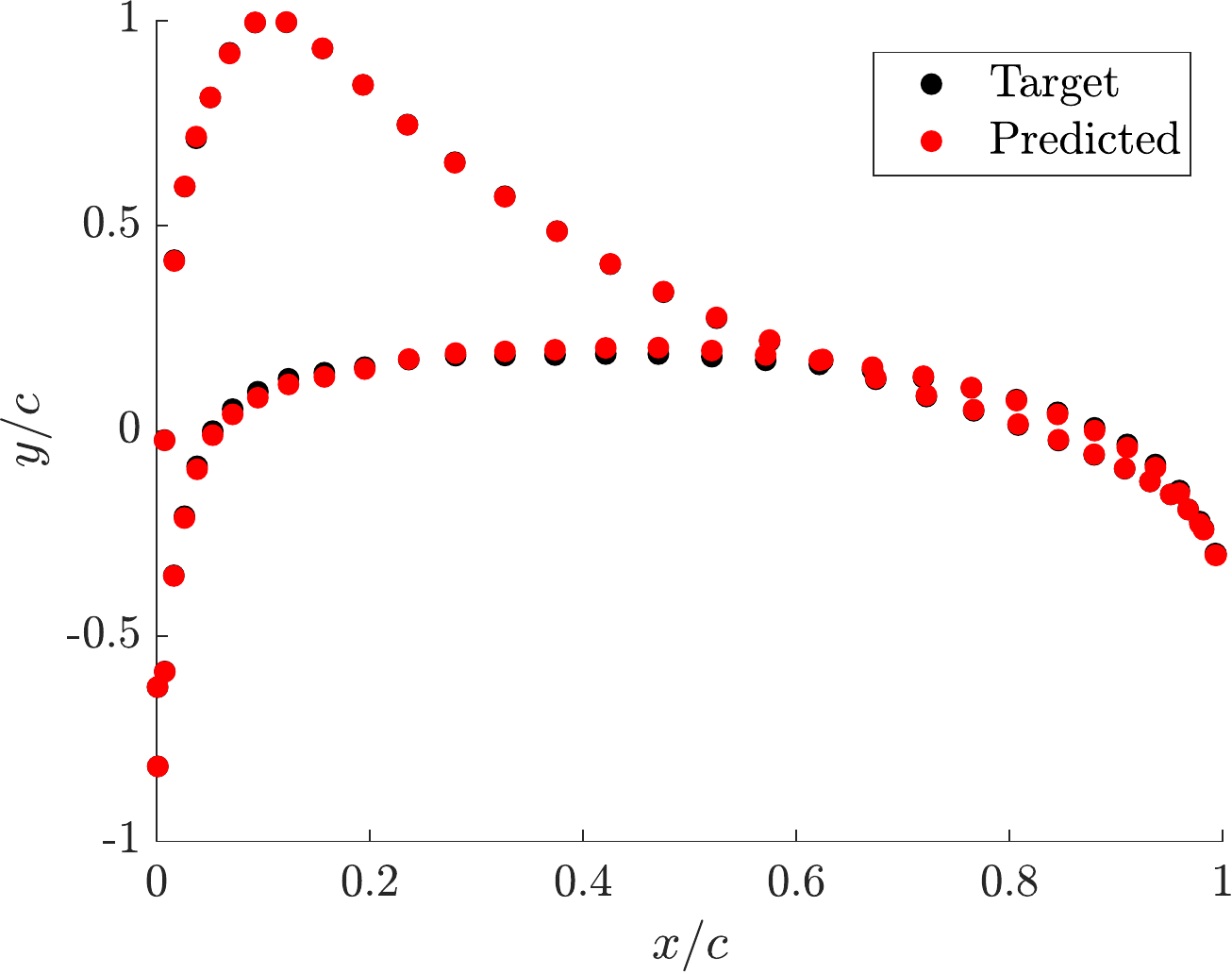}
	\caption{$C_P$ distribution comparison}
	\label{f:}
\end{subfigure} \\
\begin{subfigure}{0.5\textwidth}
	\centering
	\includegraphics[width= 1\linewidth]{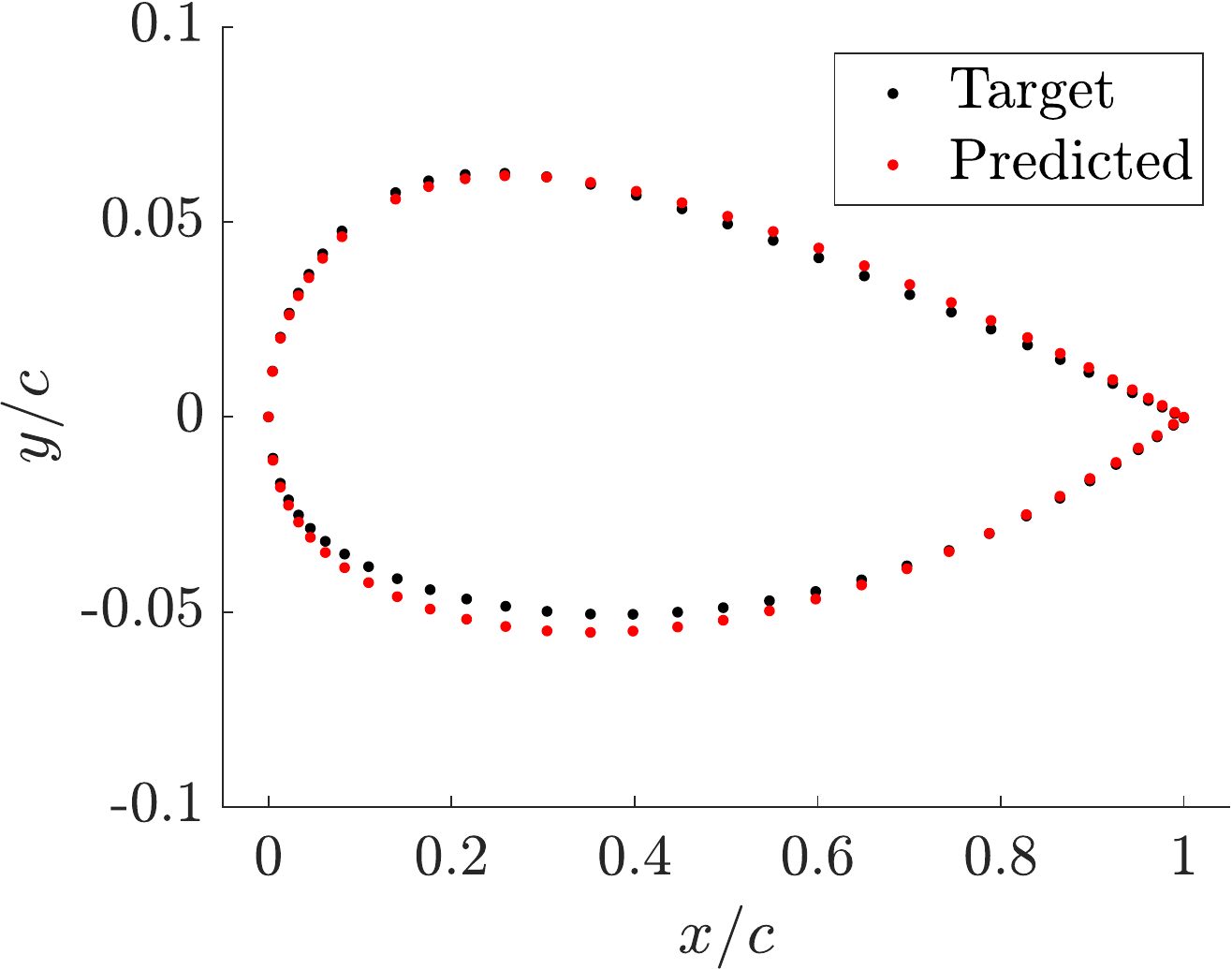}
	\caption{Airfoil shape comparison}
	\label{f:}
\end{subfigure}%
\begin{subfigure}{0.5\textwidth}
	\centering
	\includegraphics[width= 1\linewidth]{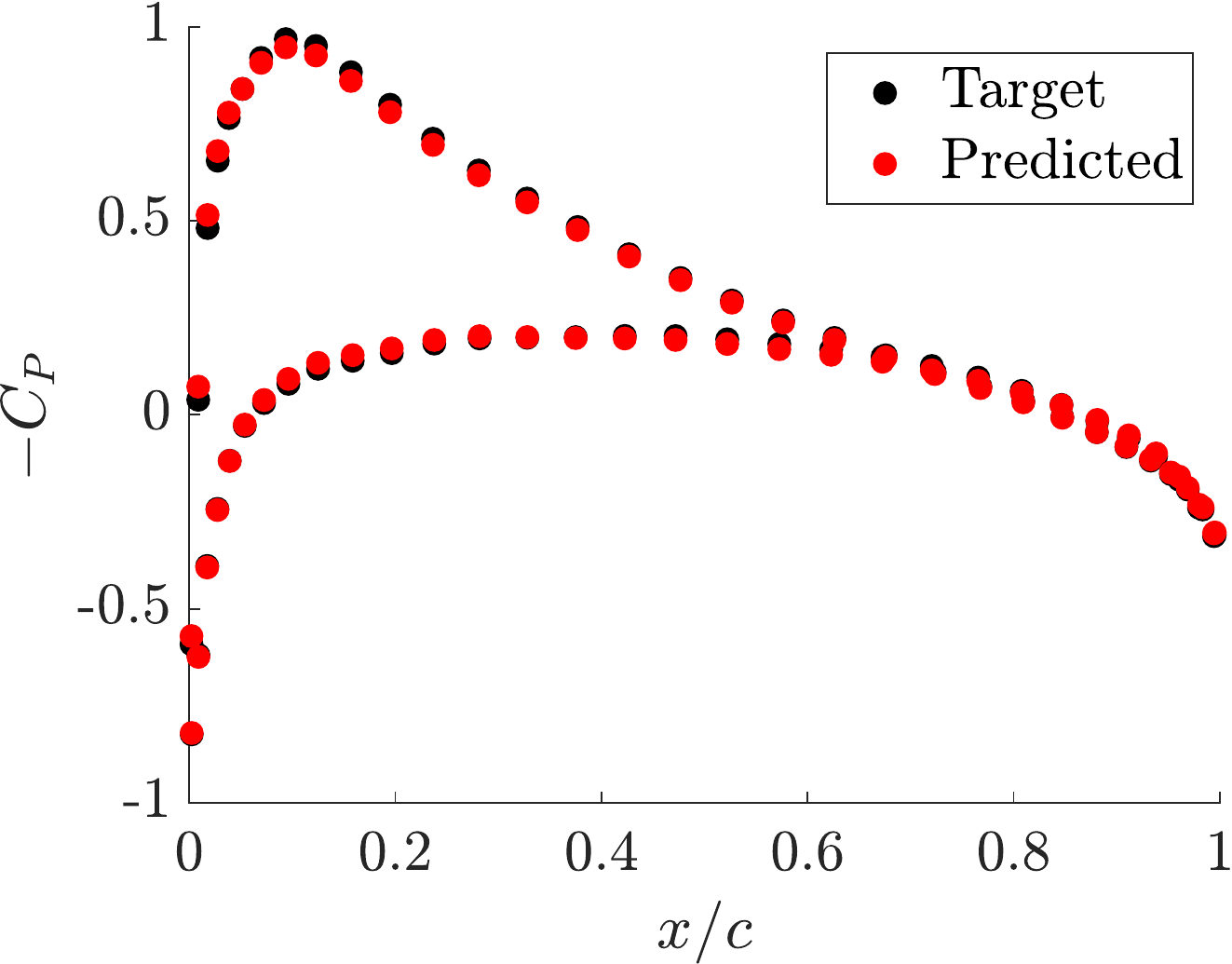}
	\caption{$C_P$ distribution comparison}
	\label{f:}
\end{subfigure}
\caption{Comparison of predicted-target design with the target. Top row: POD-Proj., bottom row: POD-Krig.}
\label{f:target_v_pred}
\end{figure}

\begin{figure}[htb!]
\centering
\includegraphics[width=.7\linewidth]{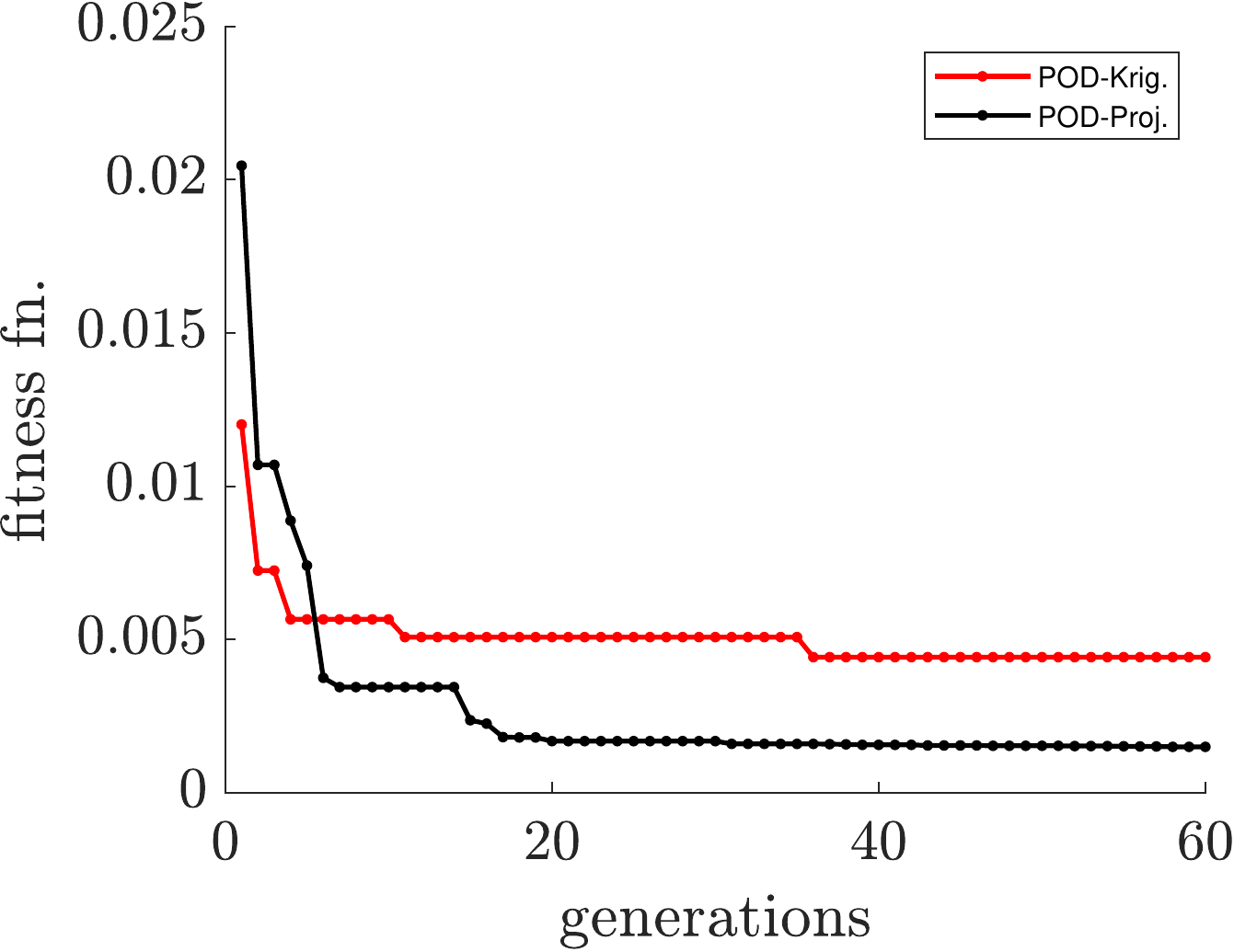}
\caption{GA optimizer convergence for inverse design problem}
\label{f:Inv_Conv_Hist}
\end{figure}

The free-stream conditions used for this test case are summarized in Table~\ref{t:Inv_des}. The optimum shape was searched using a GA based optimizer with a population size of 30 per generation and a total of 60 generations. The constraint and function convergence tolerance were set to $10^{-5}$ and $10^{-3}$ respectively and the optimization required a total of 1830 function evaluations to determine the final design. The final design is shown in Figure~\ref{f:target_v_pred}. The predictions by the proposed approach appears quite similar to the POD-Krig. However, the difference is more noticeable in the GA optimizer convergence history shown in Figure~\ref{f:Inv_Conv_Hist} where the POD-Krig. leads to a sub-optimal design compared to the proposed approach. In this specific example the proposed approach only marginally outperforms the baseline approach, but this is put to test further in predicting the lift and drag coefficients section~\ref{ss:UQ}. It was earlier shown in \cite{renganathan2018koopman} that in the transonic regime a POD-Krig. like approach is prone to predicting non-physical shock patterns mainly because the reduced state is not necessarily smoothly varying. However, such an approach is still effective if the flow does not contain parameter-dependent discontinuities as demonstrated in this example. The main benefit of the proposed approach is the gain in computational times. It required approximately 3.7 hrs of wall-clock time for the inverse design problem whereas the equivalent FOM wall-clock time for the same number of function evaluations is expected to take roughly 300 hrs.  

\subsection{Uncertainty Quantification}
\label{ss:UQ}

Next, we demonstrate the present methodology towards the uncertainty quantification problem. Uncertainty in aerospace design is inevitable and manifests itself either due to lack of knowledge (ex: biased models) or variability inherent in the system (ex: noisy manufacturing process). The aerospace design process should therefore account for the uncertainties in order to make reliable decisions early on in the design process. Specifically, we are interested in quantifying the uncertainty in the airfoil lift and drag coefficients due to the manufacturing process induced variations in the airfoil shape. We use the same $\pm~30\%$ variation in the airfoil CST coefficients and uniformly sample this design space. A Monte Carlo (MC)~\cite{mooney1997monte} simulation is carried out to propagate the input uncertainty into the model to quantify the uncertainties in the quantities of interest - namely, the lift and drag coefficients. We use the NACA test case under operating conditions defined in Table~\ref{t:Inv_des}.

\begin{figure}[H]
	\centering
	\begin{subfigure}{.5\textwidth}
		\centering
		\includegraphics[width=1\linewidth]{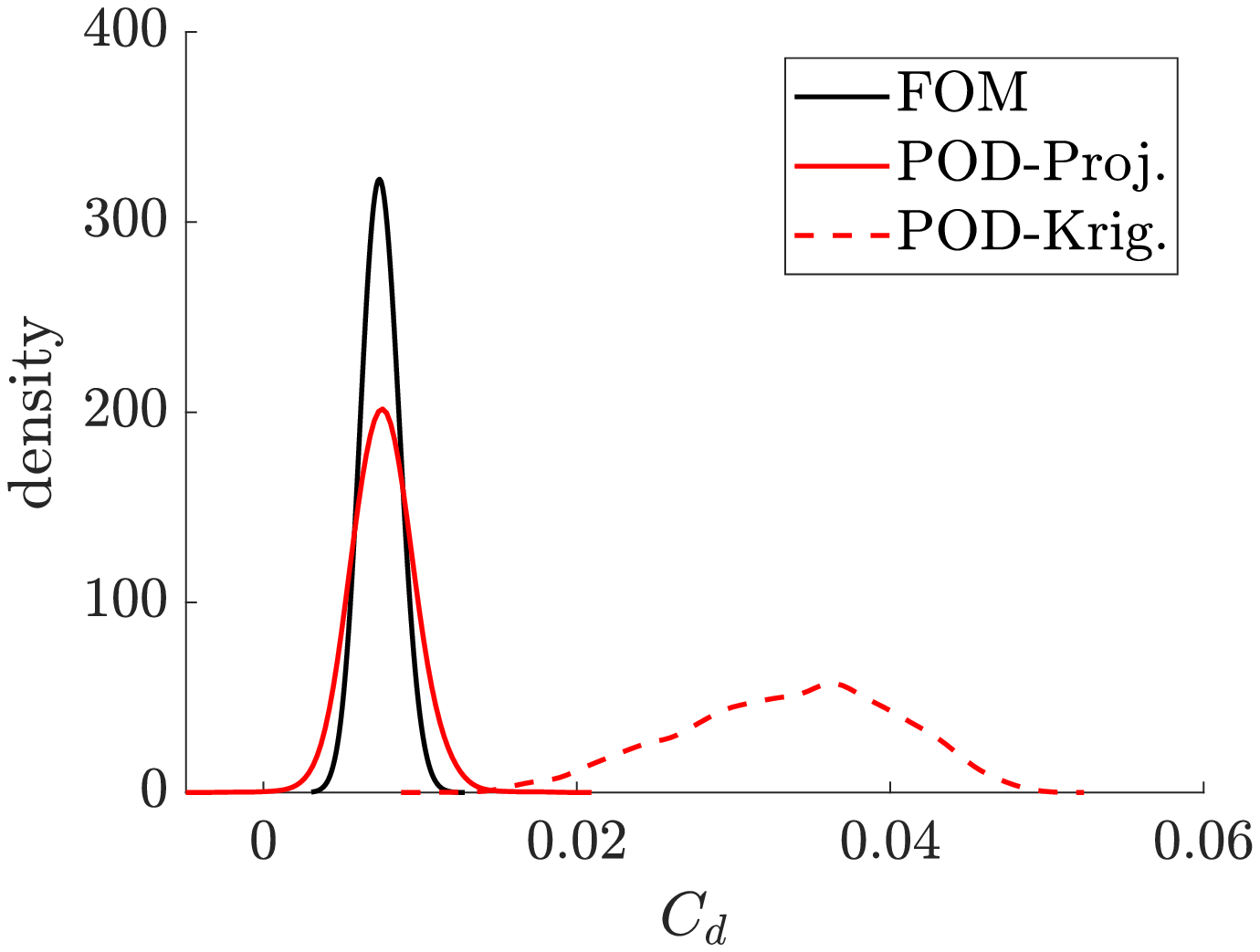}
		\caption{Drag coefficient}
	\end{subfigure}%
	\begin{subfigure}{0.5\textwidth}
		\centering
		\includegraphics[width=1\linewidth]{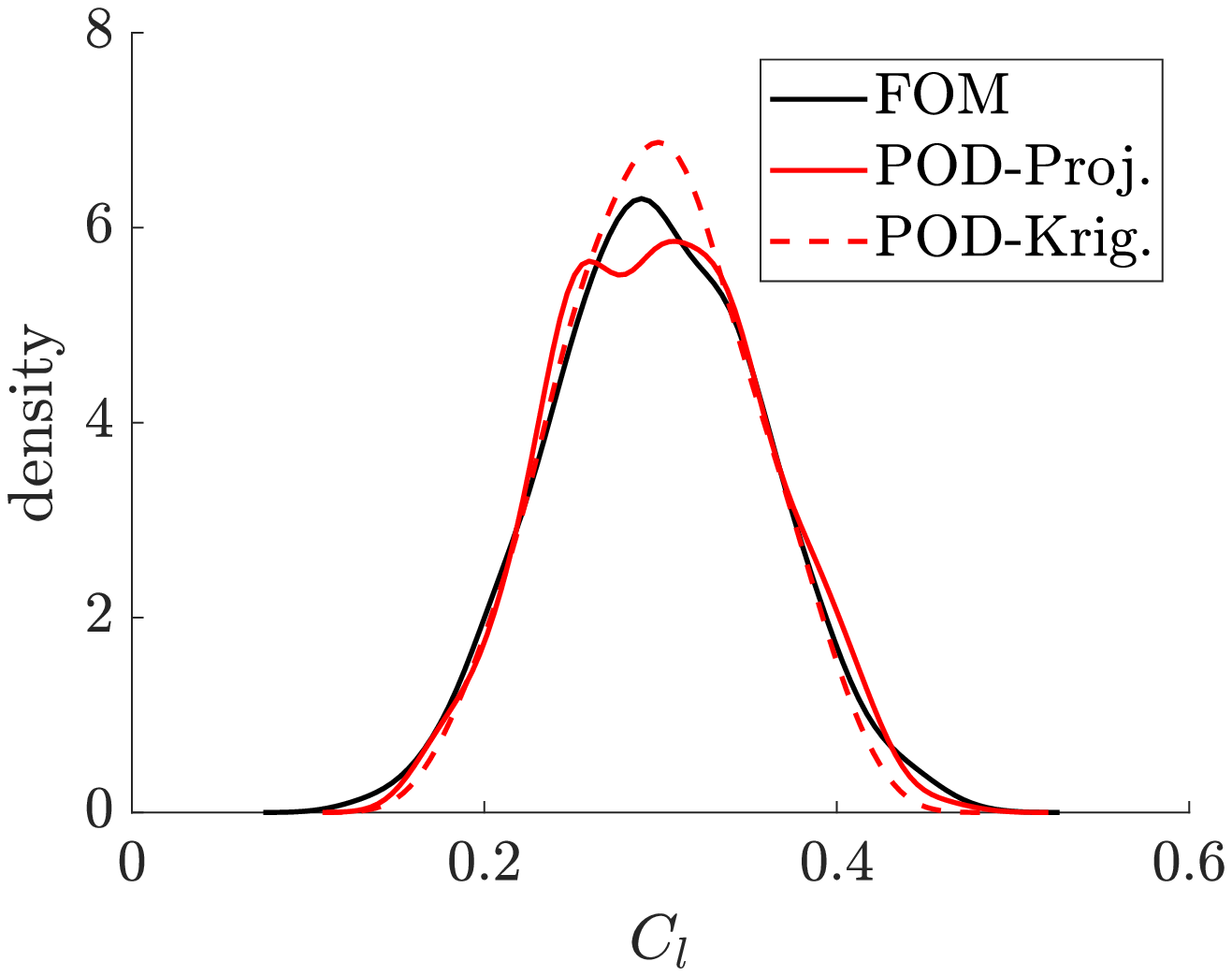}
		\caption{Lift coefficient}
	\end{subfigure}	
\caption{Kernel density estimate from 4000 monte-carlo samples.}
\label{f:kde}
\end{figure}

A total of 4000 samples were uniformly sampled from the $\pm~30\%$ range on the CST coefficients on which the proposed approach (POD-Proj.) and the baseline approach (POD-Krig.) are compared. The FOM is also evaluated at the same set of points in order to provide a reference for comparison. The probability density function is approximated from the 4000 samples via the kernel density estimation (KDE) method~\cite{epanechnikov1969non, rosenblatt1956remarks} with 100 equally-spaced points and a band-width of $0.001$ and $0.017$ for $C_d$ and $C_l$ respectively. The density plots are shown in Figure~\ref{f:kde}. In predicting the density plot for $C_l$ the two methods perform quite equally. However for $C_d$, the proposed approach clearly outperforms POD-Krig. As discussed briefly in section~\ref{A:Model_Val}, this is mainly because the $C_d$ shows greater sensitivity to $C_P$ distributions compared to $C_l$ particularly in 2-D inviscid simulations and hence even modest error in predicting $C_P$ can translate into larger errors in $C_d$. Secondly, the POD-Krig. approach does not account for the physics of the system unlike the proposed approach and hence performs poorly when the Kriging interpolator does not learn the behavior of the POD coefficients in the parameter space very well. It should be noted that there are other ways to improve upon the Kriging prediction, for instance via using a \emph{non-stationary} kernel ~\cite{paciorek2004nonstationary, rasmussen2003gaussian} to capture the underlying correlation in the data. Such models might come with a higher dimensional parametrization but have the potential to emulate rapidly varying non-linear responses very well. However, in the present work the goal is only to make a very general comparison of the proposed approach with a very basic alternative method that is applicable for non-intrusive ROM methods. The statistical measures that quantify the shape of the density plots is summarized in Table~\ref{t:clcd_stats}. As revealed by the table, the statistics for the proposed POD-Proj. approach matches the stats predicted from the FOM much better than the POD-Krig approach.

Another important aspect for discussion are the computational times required for each of the methods. All surrogate model computations were run in serial mode on a desktop computer with 4 cores and 16GB RAM. The ROM computations for 4000 samples consumed approximately 7.2 hrs of wall-clock time whereas the POD-Krig. approach consumed approximately 6.5 hrs. An equivalent budget of high-fidelity simulations run in serial for the same level of convergence as the ROM would have cost $\approx$ 667 hrs of wall-clock time; therefore the computational speed-up with the ROM is $\sim~100 \times$. However the FOM was run on a high-performance computer taking advantage of parallel computing in order to keep the wall-clock times tractable. 

 Overall, the projection-based ROM is able to capture general trends such as the range of output quantities of interest and shapes of their distributions with much better accuracy than an approach that circumvents the projection step (in this case via Kriging). This emphasizes the power of projection-based ROMs since they operate on the actual governing equations of the system as opposed to other data-driven surrogate modeling techniques. This allows them to perform more robustly with parameter variation compared to competing methods. Overall, the presented results provide evidence that the projection-based ROM has capability to potentially supplement the expensive high-fidelity models while offering compelling computational speed-ups without sacrificing the accuracy significantly.

\begin{table}[htb!]
\centering
\caption{Comparison of output statistics. Bold-face entries are significantly worse compared to rest}
\begin{tabular}{lccc|ccc}
 \hline
 & \multicolumn{3}{c}{$C_l$} & \multicolumn{3}{c}{$C_d$}\\
\hline
 & POD-Proj. & POD-Krig. & FOM & POD-Proj. & POD-Krig. & FOM \\
 \hline
Mean     &	0.2996  & 0.2968   & 0.2980  & 0.0076  & \textbf{0.0334}  & 0.0075\\
Median   &	0.2989  & 0.2967   & 0.2969  & 0.0076  & \textbf{0.0340}  & 0.0074\\
Std. Dev &  0.0587  & 0.0531   & 0.0598  & 0.0018  & \textbf{0.0068}  & 0.0007\\
Skewness &  0.0743  & \textbf{-0.0144}  & 0.0431  & 0.2294  & \textbf{-0.2765} & 0.3887\\
Kurtosis &  2.4794  & 2.4978   & 2.7697  & \textbf{4.5908}  & 2.5047  & 2.5210 \\	 
 \hline
\end{tabular}
\label{t:clcd_stats}
\end{table}

\subsection{Computational Costs}
The overall computational cost is dominated by the offline phase where the model is built. Here, we provide an estimate of the computational cost in terms of Floating Point OPerationS (FLOPS) necessary to build the ROM (POD-Proj.) as a function of the grid size ($N$) and the number of FOM snapshots ($M$). The cost of the online phase (including ROM interpolation) is trivial comparatively and the wall-clock time is more relevant in this scenario. The off-line phase includes 4 major steps whose computational cost are summarized in Table \ref{t:offline_cost};   see \cite{renganathan2018koopman} for details. Note that only the dominating factors of the cost are provided in the table. In the same table, the cost of an intrusive projection-based MOR method is also provided for comparison. It can be seen that the most expensive steps of the method are the POD and the projection which scale as $\sim~N M^2$ and $\sim~N^2 M$ respectively; as $N$ increases the cost of these steps increases quite rapidly. The finite volume discretization is relatively a cheaper step that scales linearly with grid size. Compared to a conventional intrusive MOR method, the current method incurs an additional cost due to the finite-volume discretization which can be considered a penalty paid due to the lack of access to the source code of the FOM. Whereas POD-Krig. is computationally cheaper than projection-based methods, their limitation in predicting complex flow fields have been demonstrated in this study.

\begin{table}[htb!]
\centering
\caption{Summary of offline computational cost}
\begin{tabular}{cccc}
\hline
Operation & (non-intrusive) POD-Proj. & POD-Krig. & (intrusive) POD-Proj.\\
\hline
\hline
Snapshot Scaling & $\mathcal{O}(N)$ & $\mathcal{O}(N)$ & $\mathcal{O}(N)$\\
POD & $\mathcal{O}(NM^2)$ & $\mathcal{O}(NM^2)$ & $\mathcal{O}(NM^2)$\\
Finite Vol. Discret. & $\mathcal{O}(MN)$ & - & - \\
Projection & $\mathcal{O}(N^2 M)$ & - & $\mathcal{O}(N^2 M)$\\
\hline
\end{tabular}
\label{t:offline_cost}
\end{table}
The wall-clock times for online evaluation are summarized in Table~\ref{t:wall-clocK}, where the computational gains of the proposed approach are compelling  and therefore given their demonstrated accuracy in predicting non-linear flows, are viable candidates during early and middle stages of aerospace design. While the POD-Krig. approach requires a comparable amount of wall-clock time as the proposed approach, they are based on prediction for only one variable ($p$) whereas the proposed approach solves for $O=8$ observables. Therefore for a fair comparison, the POD-Krig. wall-clock times are multiplied by $O$ in Table~\ref{t:wall-clocK} to provide the perspective. Having said that, when surrogates of only one field variable (such as pressure) is of interest, this study suggests that a simpler method such as the POD-Krig. could be sufficient. 

\begin{table}[htb]
    \centering
        \caption{Comparison of computational (wall-clock) times}
    \begin{tabular}{cccc}
        \hline\hline
         & FOM & POD-Krig. & POD-Proj. \\
        \hline
        Evaluation        &  600s & 5.85s$(\times O)$& 6.48s \\
        Inverse Design    &  300hrs & 3.4hrs$(\times O)$ & 3.7hrs \\ 
        Uncert. Quant.&  667hrs & 6.5hrs$(\times O)$ & 7.2hrs \\
        \hline
    \end{tabular}
    \label{t:wall-clocK}
\end{table}

\section{Conclusion}

We extend the earlier work by ~\cite{renganathan2018koopman} to systems with parametric geometry in addition to parametric boundary conditions, which are more commonly encountered in the aerospace design process. The methodology is validated under subsonic and transonic compressible inviscid flow. Under subsonic conditions, exceptional accuracy (<5\%) was observed in approximating the flow-field as well as ouputs $C_l$, $C_d$ and $C_P$, whereas in the transonic regime the average error is higher (in the 15-20\% range). The method suffers from typical limitations of a POD-based methods to capture highly non-linear flow such as ones with moving shocks~\cite{Beran2001, legresley2006application}. A common trend observed in this work is that $C_P$ and $C_l$ are predicted with much better accuracy than $C_d$, albeit showing better accuracy and robustness than a competing non-intrusive method.

Furthermore to demonstrate its computational efficiency, the method is applied to 2 specific applications in the many-query context: (i) inverse design and (ii) uncertainty quantification via Monte Carlo sampling. In both situations, the ROM is executed $\mathcal{O}(1000)$ times at a wall-clock time of 2-8 hrs, while the equivalent budget of FOM would have consumed 200-800 hrs. Therefore first and foremost, these results demonstrate the utility of the methodology towards real-time decision making.  The inverse design test case is mainly the test of the ROM to satisfy the physics of the problem, even in its approximated form. The predictions of the ROM is within $5\%$ of the target solution while again, achieving it at superior computational efficiency compared to the FOM. Lastly, a Monte Carlo analysis with 4000 uniformly sampled points from the input space was used to approximate the probability distributions of the two main outputs considered in this work: the $C_d$ and $C_l$. Overall, the $C_l$ showed better match with the FOM results, while the $C_d$ predicted the statistics with relatively higher discrepancy. Regardless, the predictions of the ROM turn out to capture general trends at a fraction of the computational cost of the FOM, while again showing more robustness than a competing non-intrusive method.

The present work is a first step towards performing projection-based model reduction with black-box models. Therefore to keep the exposition simple, the 2-dimensional compressible euler equations are used as the test case. However, the methodology naturally extends to 3-dimensional systems with $S>4$ without modification. The only requirement is complete knowledge of the governing equations in continuous PDE form, which is always available in the form of code documentation with black-box commercial codes. The author is currently investigating purely data-driven methods where this requirement can be waived. Another direction for future work is the adaptive construction of the model, where snapshots are sequentially generated based on certain \emph{greedy} criterion. Finally, state-of-the art methods to enhance the performance of the method for flows with discontinuities is also being investigated.

Overall, the present methodology establishes feasibility for projection based MOR for static parametric systems available as a black-box in addition to demonstrating their utility towards aerospace design. Comparison against the best known non-intrusive MOR technique (given the same constraints as the present study), namely POD+Interpolation suggested that the projection-based approach leads to more accurate predictions in the presence of parameter-dependent discontinuities. However, the author acknowledges that projection-based ROMs, in addition to incurring a greater off-line cost, can also lead to numerical stability issues (specifically in time-dependent systems) which needs to be addressed. Therefore the choice between either of these methods is dependent upon the flow regime under consideration and hence the domain knowledge of the engineers has to be leveraged.

\clearpage
\section*{Acknowledgments}

\bibliography{main}

\section{Appendix}
\subsection{Shape parameterization}
\label{A:Shape}

The CST model of parameterization defines a \textit{class} function $C$ and a \textit{shape} function $S$ and the curve being parameterized is specified as their product. The class function, $C(\psi)$ is more generally defined as

\begin{equation*}
C_{n_1}^{n_2} (\psi) := \psi ^{n_1}(1 - \psi)^{n_2}
\label{e:class_fn}
\end{equation*}

where the variable $\psi$ represents the non-dimensional chord-wise distance, $n_1$ and $n_2$ define the specific class. The unique shape of an airfoil is driven by the shape function, specified as follows

\begin{equation}
S(\psi) = \sum_{i=0}^{n} A_i \psi^i 
\label{e:shape_fn}
\end{equation}

It is particularly useful to define a \textit{unit shape function}, i.e. $S(\psi) = 1$ such that the individual coefficients $A_i$ \footnote{ $A_i$'s are denoted by $K_{i,n}$ for unit shape functions} can be obtained as generic constants. For instance for $n=1$ the simplest decomposition one could get for the shape function is $S(\psi) = S_0(\psi) + S_1(\psi)$ where $S_1(\psi) = \psi$ and $S_2(\psi) = 1- \psi$ where the coefficients $A_0=1$ and $A_1 =1$. Similarly, for the general $nth$ order shape function, the decomposition of the unit shape function can be done using $\textit{Bernstein}$ polynomials 

\begin{equation}
S(\psi) = \sum_{i=0}^n K_{i,n} \psi^i (1 - \psi)^{n-i}
\label{e:Bernstein}
\end{equation}

where the coefficients are the \textit{binomial} coefficients given by

\[ K_{i, n} = {n \choose i} = \f{n!}{i!~(n-i)!}\]

The final shape of the airfoil shape is then given by

\begin{equation}
\mathbf{y}(\psi) = C(\psi)S(\psi)
\label{e:CST_final}
\end{equation}

\begin{figure}[htb!]
\centering
\includegraphics[width=1\textwidth]{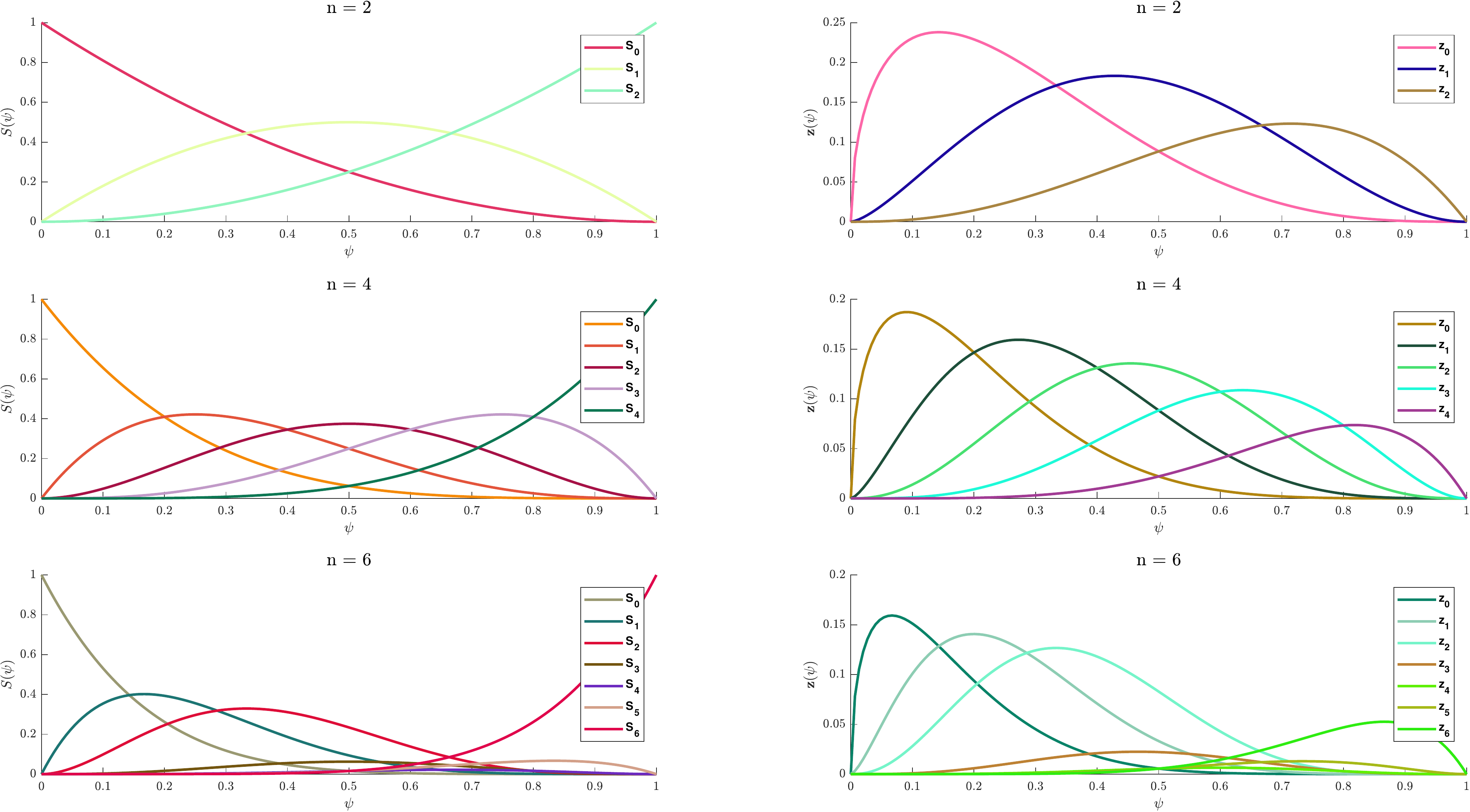}
\caption{Examples of the shape function decomposition (into Bernstein polynomials) for various values of the order $n$ (\textit{left}) and the resulting component airfoils (\textit{right}). The coefficients $A_i$ correspond to unit shape function, which can be scaled up/down to obtain a specific airfoil shape.}
\label{f:ShapeFnDemo}
\end{figure}

The unit shape functions and the corresponding airfoil geometries are illustrated in the Figure~\ref{f:ShapeFnDemo}. It can be seen that such a parametrization results in each component shape function peak being equally distributed between the leading and trailing edges leading to the same effect in the component airfoils. It is now a matter of scaling up or down, the binomial coefficients of the Bernstein polynomials in order to approximate the unique airfoil shape of interest.

The coefficients $A_i$ represent the actual parameters of the shape, given $n$ the order of the Bernstein polynomials. An $nth$ order CST parametrization has $n+1$ parameters. If separate parametrizations are sought for the upper and lower surfaces of the airfoil, then the CST parametrization leads to $2(n+1)$ parameters to specify the whole shape of the airfoil, where the $n$ needs to be determined for a specific geometry under consideration. However, typically $n=3-5$ are observed to be adequate to parametrize the airfoil shapes considered in this work. One way to determine $n$ and the associated polynomial coefficients is to find the values that minimizes certain error between the true shape of the airfoil and the resulting approximation via CST. In this work, the parameters for a given airfoil shape are determined by solving the following minimization problem

\begin{equation}
\underbrace{\text{minimize}}_{A_j, ~n}~ \left \| \mbf{y}(\tilde{\psi}) - C(\tilde{\psi})S(\tilde{\psi},~A_j) \right \|_2 ^2 
\label{e:CST_param}
\end{equation}

where $\tilde{\psi} \in \R^{n+1}$ are $n+1$ equally spaced points sampled from $\psi$ spanning $[0,1]$ \footnote{Note that picking equally spaced points works well for the current airfoil geometries under consideration in this work and may not work for any arbitrary geometry. More generically, a least-squares fit considering all the points might be more suitable}. This way, the smallest possible $n$ and their corresponding Bernstein coefficients are determined. For the RAE2822 airfoil shape, the following parameterization was obtained ($n=3$): 

\[ A_{RAE2822} = \begin{bmatrix}
 0.1268 & 0.4670 & 0.5834 & 0.2103 \\
-0.1268 &-0.5425 &-0.5096 & 0.0581
\end{bmatrix} \]

 where the first and second rows represent the parameterization of the upper \& lower surfaces of the airfoil; the comparison of the CST curve and the actual RAE2822 shape is shown in Figure~\ref{f:RAE2822_CST}. It can be seen that the CST parametrization gives an \textit{adequate} approximation to the true curve with only 8 parameters. The coefficients may now be perturbed to modify the baseline airfoil shape.

 Similarly, the NACA0012 airfoil shape is approximated via CST and is also shown in Figure~\ref{f:NACA0012_CST}. In this case, due to the lack of camber, the CST gives very good approximation with $n=2$. Additionally, due to the symmetry of the airfoil about the chord, the parameters (given below) are equal in magnitude and opposite in sign. However, all the $2(n+1)$ degrees of freedom are considered in this work for the NACA0012 airfoil.

\[ A_{NACA0012} = \begin{bmatrix}
 0.1689 & 0.2699 & 0.1387 \\
-0.1689 &-0.2699 &-0.1387
\end{bmatrix} \]
\newpage

\subsection{Model validation data}
\label{A:Model_Val}

The values of the output quantities of interest and their associated errors (defined by Eq.~\ref{e:error_metrics}) is presented in Tables~\ref{t:NACA_Shape_CLCD} and \ref{t:RAE_Shape_CLCD} below.

\begin{table}
\centering
\caption{Comparison of $C_P$, $C_l$ and $C_d$ between ROM \& FOM for the NACA0012 test case}
\begin{tabular}{cc|ccc}
\hline
\hline
Case & $C_P$ Error \% & $C_l$ (ROM) & $C_l$ (FOM) & Error \% \\
\hline
1  & 1.29 & 0.1889 & 0.1912 & 1.20  \\
2  & 0.74 & 0.2018 & 0.2070 & 2.50  \\
3  & 0.80 & 0.2932 & 0.2943 & 0.37  \\
4  & 1.86 & 0.2795 & 0.2865 & 2.44  \\
5  & 1.36 & 0.3550 & 0.3621 & 1.96  \\
6  & 0.62 & 0.3691 & 0.3664 & 0.73  \\
7  & 0.46 & 0.3298 & 0.3272 & 0.79  \\
8  & 2.79 & 0.3229 & 0.3312 & 2.50  \\
9  & 0.76 & 0.3109 & 0.3137 & 0.89  \\
10 & 5.43 & 0.2710 & 0.3065 &11.58  \\
\hline
\end{tabular}
\label{t:NACA_Shape_CLCD}
\end{table}

\begin{table}
\centering
\caption{Comparison of $C_P$, $C_l$ and $C_d$ between ROM \& FOM for the RAE2822 test case}
\begin{tabular}{cc|ccc|ccc}
\hline
\hline
Case & $C_P$ Error \% & $C_d$ (ROM) & $C_d$ (FOM) & Error \% & $C_l$ (ROM) & $C_l$ (FOM) & Error \% \\
\hline 
1	&	7.68	&	0.0161	&	0.0174	&	7.47	&	0.9174	&	0.9608	&	4.52	\\
2	&	12.93	&	0.0336	&	0.0302	&	11.26	&	1.0669	&	0.9825	&   8.59	\\
3	&	2.84	&	0.0257	&	0.0262	&	1.91	&	1.1140	&	1.1446	&	2.67	\\
4	&	8.14	&	0.0321	&	0.0264	&	21.59	&	1.1004	&	1.0901	&	0.95	\\
5	&	13.77	&	0.0194	&	0.0224	&	13.39	&	0.8965	&	1.0345	&	13.34	\\
6	&	12.78	&	0.0298	&	0.0484	&	38.43	&	1.0288	&	1.0082	&	2.04	\\
7	&	5.91	&	0.0306	&	0.0279	&	 9.68	&	0.9403	&	0.9286	&	1.26	\\
8	&	4.99	&	0.0245	&	0.0245	&	 0	&	0.9124	&	0.8914	&	2.36	\\
9	&	5.63	&	0.0292	&	0.0315	&	7.30	&	0.9507	&	0.9779	&	2.78	\\
10	&	9.96	&	0.0326	&	0.0217	&	50.23	&	0.9750	&	0.9793	&   0.44	\\
\hline
\end{tabular}
\label{t:RAE_Shape_CLCD}
\end{table}
\clearpage

\subsection{Discrete Empirical Interpolation Method (DEIM)}
\label{A:DEIM}

The Discrete Empirical Interpolation Method (DEIM) is briefly reviewed here and as an illustration one of the non-linear constraints used in Eq.~\ref{e:Euler_Cons} is evaluated. For a non-linear function $\mbf{f}(\theta)\in \R^N$ the DEIM approximates $\mbf{f}$ by projecting it onto a subspace spanned by $\lbrace \mbf{x}_1,...,\mbf{x}_q\rbrace\subset \R^N$ as
\begin{equation}
\mbf{f}(\theta) \approx \mbf{X} c(\theta)
\end{equation}

where $\mbf{X}=[\mbf{x_1},...,\mbf{x_q}] \in \R^{N \times q},~~q<<N$ is determined via a POD of the snapshots of $\mbf{f}$ and is assumed to be globally valid in the design space that bounds the design parameters $\theta$ and $\mbf{c}(\theta)\in\R^q$ are the coefficients of the basis expansion. Then the approximation of $\mbf{f}$ requires only the determination of $\mbf{c}(\theta)$ which requires only $q$ equations. The DEIM gives a distinguished set of $q$ points from the over-determined system $\mbf{f}(\theta) = \mbf{X}\mbf{c}(\theta)$. Given a permutation matrix $\mbf{P}$ that would give $q$ such distinguished rows of a matrix when pre-multiplied, then the $q\times q$ system necessary to solve for the coefficients is given by

\begin{equation}
\mbf{P}^\top \mbf{f}(\theta) = (\mbf{P}^\top\mbf{X})\mbf{c}(\theta)
\end{equation}

So the approximation of $\mbf{f}(\theta)$ is then given by

\begin{equation}
\mbf{f}(\theta) \approx \mbf{X}(\mbf{P}^\top\mbf{X})^{-1} \mbf{P}^\top \mbf{f}(\theta)
\end{equation}

If the $q$ row-indices (that are extracted by pre-multiplying with $\mbf{P}^\top$) are represented by a vector, $\mbf{\varrho}$, then in the above equation, $\mbf{P}^\top \mbf{f}(\theta)$ is equivalent to extracting the $\varrho$ rows of $\mbf{f}$. Therefore the approximation of $\mbf{f}(\theta)$ requires only $q$ computations which is efficient because $q<<N$. Similarly, a non-linear function that depends on the state, $\mbf{f}(\mbf{u})$ can be approximated as

\begin{equation}
\mbf{f}(\mbf{u}) \approx \mbf{X} (\mbf{P}^\top\mbf{X})^{-1} \mbf{P}^\top \mbf{f}(\mbf{u})
\end{equation}

Since $\mbf{u} = \Phi_k^\top \tilde{\mbf{u}}$ and setting $\tilde{\mbf{f}}=\Phi_k^\top \mbf{f}(\mbf{u})$, $\tilde{\mbf{f}}$ can be approximated as

\begin{equation}
\tilde{\mbf{f}} = \Phi_k^\top \mbf{X} (\mbf{P^\top}\mbf{X})^{-1} \mbf{f}(\mbf{P}^\top \Phi_k \tilde{\mbf{u}})
\end{equation}

In the above equation, the term $\Phi_k^\top \mbf{X} (\mbf{P^\top}\mbf{X})^{-1}$ is independent of the state and hence can be pre-computed and $\mbf{P}^\top \Phi_k$ is just extraction of the $\varrho$ rows of $\Phi_k$. Therefore using the DEIM, the non-linear term can be expressed in terms of the reduced state, $\tilde{\mbf{u}}$ and hence can be efficiently computed.

Now the DEIM is illustrated on evaluating the first constraint of Equation~\ref{e:Euler_Cons} which in discretized form is given below

\begin{equation}
\mbf{h}_1 = \mbf{y}_5 - \frac{\mbf{y}_1\mbf{y}_3}{\mbf{y}_2}
\end{equation}

Let $\varrho_5$ be the vector containing the $q$ row-indices returned by DEIM via snapshots of the non-linear term $\mbf{y}_5$ and $\Phi_1$, $\Phi_2$, $\Phi_3$, $\Phi_5$ be the projection matrix of $\mbf{y}_1$, $\mbf{y}_2$, $\mbf{y}_3$ and $\mbf{y}_5$ respectively. Then

\begin{equation}
\tilde{\mbf{h}}_1 = \tilde{\mbf{y}_5} - \Phi_5^\top \mbf{X} ~[\mbf{X}(\varrho_5,:)]^{-1} \left\lbrace \frac{\Phi_1(\varrho_5,:)\tilde{\mbf{y}_1}~~\Phi_3(\varrho_5,:)\tilde{\mbf{y}_3}}{\Phi_2(\varrho_5,:)\tilde{\mbf{y}_2}} \right\rbrace
\end{equation}

In the above equation, the term outside of the braces can be pre-computed. Additionally since $\mbf{y}_5 = \frac{\mbf{y}_1\mbf{y}_3}{\mbf{y}_2}$, $\mbf{X}=\Phi_5$ and hence the term reduces to $[\mbf{X}(\varrho_5,:)]^{-1}$ which is $q \times q$ and hence can be cheaply computed. Therefore using the DEIM, the non-linear constraints are evaluated in terms of the reduced state variables which makes it computationally cheap.

\end{document}